\documentclass[11pt]{article} 
\usepackage{amsmath, amscd, amssymb,
  amsthm, verbatim, amsrefs}
\usepackage{url}
\usepackage[ocgcolorlinks,linktoc=all]{hyperref}
\hypersetup{citecolor=blue,linkcolor=red}
\allowdisplaybreaks[4]

\newtheorem{theorem}{Theorem}[section]

\newtheorem{lemma}[theorem]{Lemma}
\newtheorem{proposition}[theorem]{Proposition}

\newtheorem{definition}[theorem]{Definition}
\newtheorem{remark}[theorem]{Remark}

\newtheorem{assumption}[theorem]{Assumption}

\newcommand{\T}{\mathbb{T}}
\newcommand{\n}{\nabla}
\newcommand{\p}{\mathbb{P}}
\newcommand{\no}{\mathbf{n}}
\newcommand{\D}{\mathrm{div}}
\newcommand{\Q}{\mathbf{Q}(\n u)}
\newcommand{\V}{\mathbf{v}(\n u)}

\newcommand{\E}{\mathbb{E}} 

\newcommand{\N}{\mathbb{N}} 
\newcommand{\R}{\mathbb{R}} 
\newcommand{\Div}{\mathrm{div}}

\newcommand\newdot{{\kern.8pt\cdot\kern.8pt}}

\begin{document}

\title{Existence of martingale solutions  for a stochastic weighted mean curvature flow of graphs}

\author{Qi Yan \footnote{Academy of Mathematics
and Systems Science, Chinese Academy of Sciences, No. 55, Zhongguancun East Road, Beijing,
100190.  \texttt{yanqi19@mails.ucas.ac.cn}} \quad and  \quad Xiang-Dong Li \footnote{1. State Key Laboratory of Mathematical Sciences, Academy of Mathematics
and Systems Science, Chinese Academy of Sciences, No. 55, Zhongguancun East Road, Beijing,
100190, China; 2. School of Mathematical Sciences, University of Chinese Academy of Sciences,
Beijing, 100049, China. Research supported by National Key R\&D Program of China (No. 2020YF0712702), NSFC No.
 12171458, and Key Laboratory RCSDS, CAS, No. 2008DP173182.  \texttt{xdli@amt.ac.cn}}    }
  
\date{\today}

\renewcommand{\theequation}{\arabic{section}.\arabic{equation}}
\maketitle

\begin{abstract}
  \noindent
  We are concerned with a stochastic mean curvature flow of graphs with extra force over a periodic domain of any dimension. Based on compact embedding method of variational SPDE, we prove the existence of martingale solution. Moreover, we derive the  small perturbation limit of the stochastic weighted mean curvature flow.  
\end{abstract}

\section{Introduction}

The mean curvature flow(MCF) is a significant geometric flow of hypersurfaces both in theory and application, which describes the growth of crystal, the evolution of interface of two phases and so on. 

Given a family of smooth $n$-dimensional hypersurfaces $M_t$ in $\R^{n+1}$, the mean curvature flow is characterized by 
\[
  \partial_t x=-H\no 
\]
where $x$ is the position vector, $\no$ is the outer normal at $x\in M_t$ and $H$ denotes the mean curvature of the hypersurface at $x\in M_t$, see \cite{Huisken1984}.

MCF is the gradient flow of area functional in geometry. Suppose that $M$ is a closed hypersurface in $\R^{n+1}$ and $M_t$ is a variation of $M$. That is, $M_t$ is a one-parameter family of hypersurfaces with $M_0 = M$. If we think of volume as a functional on the space of hypersurfaces, then the ﬁrst variation formula gives the derivative of volume under the variation
\[
  \frac{d}{d t} \operatorname{Area}\left(M_t\right)=\int_{M_t}\left\langle\partial_t x, H \mathbf{n}\right\rangle.
\]
It follows that the gradient of volume is
\[
  \nabla \mathrm{Area}=H \mathbf{n}.
\]
Thus, we have 
\[
  \partial_t x=-H\no=-\n \mathrm{Area},
\]
which means that MCF is the gradient flow of area functional (with $L^2$-metric), see \cite{Colding2010}.

The weighted area with respect to the weighted measure $d\mu=e^{-f}dx$ is 
\[
  \mathrm{Area}_f(M_t)=\int_{M_t}e^{-f}dx,
\]
where $f$ is a smooth function on $\R^{n+1}$. By the first variation formula, we have
\[
  \frac{d}{d t} \operatorname{Area}_f\left(M_t\right)=\int_{M_t}\left\langle\partial_t x, H_f \mathbf{n}\right\rangle d\mu,
\]
where $H_f=H-\langle \overline{\n} f,\no\rangle$ is called the $f$-mean curvature, and $\overline{\nabla}$ denotes the gradient in $\R^{n+1}$. The gradient flow of the weighted volume:
\[
  \partial_t x=-H_f\no.
\]
is called  the \textbf{$f$-mean curvature flow}, see \cites{Wei09, Jian2007, JL08, Jian2009, Liu2012, LW2015, meira_space_2020} for more information about $f$-mean curvature flow. When $f=\frac{C|x|^2}{2}$ (where $C>0$), this flow is called the Gauss mean curvature flow; see \cites{BM10,BM14}. 
When $f$ is constant, the $f$-mean curvature flow reduces to the mean curvature flow. In this setting, the $f$-mean curvature flow can be viewed as the mean curvature flow on the smooth metric measure space $(\R^{n+1},dx^2,d\mu)$. The $f$-mean curvature flow can also be interpreted as a mean curvature flow. Let $F_t:M^n\to\R^{n+1}$ be an $f$-mean curvature flow, and let $dx^2$ be the Euclidean metric on $\R^{n+1}$. Consider the mapping:
\begin{equation*}\label{flow2}
	\begin{array}[pos]{rcl}
		\widetilde{F}_t:M^n\times\R&\to&\R^{n+1}\times\R\\
		(x,s)&\mapsto&(F_t(x),s).
	\end{array}
\end{equation*}
where the metric on $\R^{n+1}\times\R$ is $dx^2+e^{2f(x)}ds^2.$ Let $\widetilde{M}_t=M_t\times\R$. According to \cite{Smoczyk01}*{Corollary 2.4}, we have
\[
	\mathbf{H}_{\widetilde{M}_t}(x,s)=\left(\mathbf{H}_f(x),0\right),
\]
where $\mathbf{H}_f=-H_f\mathbf{n}$ is the $f$-mean curvature vector of $M_t=F_t(M)$, and $\mathbf{H}_{\widetilde{M}_t}$ is the mean curvature vector of $\widetilde{M}_t$ in the Riemannian manifold $(\mathbb{R}^{n+1}\times\R,dx^2+e^{2f(x)}ds^2)$. Consequently, we obtain
\[
	\frac{\partial\widetilde{F}}{\partial t}(x,s)=\left(\frac{\partial F}{\partial t},0\right)=\left(\mathbf{H}_f,0 \right)=\mathbf{H}_{\widetilde{M}_t}(x,s).
\]
Thus, the $f$-mean curvature flow is the mean curvature flow in the Riemannian manifold $(\mathbb{R}^{n+1}\times\R,dx^2+e^{2f(x)}ds^2)$.
Furthermore, the $f$-mean curvature flow also arises in the study of Ginzburg-Landau vortex dynamics:
\[
\frac{\partial V_{\epsilon}}{\partial t}=\overline{\Delta}V_{\epsilon}+\overline{\nabla}f\overline{\nabla}V_{\epsilon}+AV_{\epsilon}+\frac{BV_{\epsilon}}{\epsilon^2}(1-|V_{\epsilon}|^2),\quad V_{\epsilon}:\mathbb{R}^n\times [0,T]\to\mathbb{R}^2.
\]
For details, see \cite{JL06}. 

Mean curvature flow with a random force was first proposed by Kawasaki and Ohta in \cite{kawasaki1982} from time-dependent Ginzburg-Landau model incorporating the influence of thermal noise. Yip \cite{Yip1998} obtained the existence of subsets whose boundaries evolve in the normal direction with velocity formally equal to the mean curvature plus a random driving force. The randomness is introduced by means of stochastic ﬂows of diffeomorphisms generated by Brownian vector ﬁelds which are white in time but smooth in space. The result was obtained by means of a time-stepping scheme, in which in one step the set evolves by pure mean curvature, and the evolution is defined through a variational minimization, and in the next step the set is transported by a random flow of diffeomorphisms. Yip finally proved the tightness of the laws of approximation of subsets. Dirr, Luckhaus and Novaga \cite{Dirr2001} , Souganidis and Yip \cite{Souganidis2004} independently studied  the stochastic mean curvature flow of two tangential unit circles on the plane with a small random force and proved a selection principle that in the limit of vanishing the intensity of the random force the measure can select two solutions, the upper and lower barrier in the sense of De Giorgi respectively, viscosity supersolution and subsolution in the sense of Chen-Giga-Goto \cite{Giga}. Es-Sarhir and von Renesse \cite{Renesse2012} studied the  stochastic curve shortening flow of graphs, which is stochastic MCF for curves in the plane, and obtained the existence and uniqueness of strong solution and its long-time behavior by the classical framework of SPDEs with the standard coercivity assumption replaced by a so-called Lyapunov condition. Hofmanov\'{a}, R\"{o}ger and von Renesse \cite{Hofmanova2017} proved the existence of weak martingale solution of stochastic mean curvature flow of two-dimensional graphs using three-step approximation scheme (that is, adding a small viscous term to make the equation uniformly parabolic, adding a small-coefficient high enough order term to yield a semilinear nondegenerate parabolic SPDE and truncating the non-Lipschitz continuous term to use the mild solution method), fixed point method and the Gauss-Bonnet formula of two-dimensional torus. Dabrock, Hofmanov\'{a} and R\"{o}ger \cite{Dabrock2021} obatined the existence of martingale solution, furthermore, its longtime behavior which extended the result of \cite{Renesse2012}*{Theorem 4.2}, of stochastic MCF of graphs of all dimensions using the compact embedding method in variational SPDEs see \cites{Viot, Pardoux2021, daprato_zabczyk_2014, Gawarecki2011}. A uniform energy estimate plays a significant role in their proof of the tightness of viscous approximation.

In this paper, we study the so-called $f$-mean curvature flow of $n$-dimensional hypersurfaces given as graphs over $n$-dimensional torus $\T^n$ with a stochastic perturbation. We consider a random evolution of hypersurfaces $(M_t)_{t\geq 0}$ in $\R^{n+1}$ given by immersions $\phi_t:M\to\R^{n+1}$, where $M$ is a smooth manifold, and a real-valued Wiener process $W$ on the probability  space $(\Omega,\mathcal{F},\p)$. $\phi_t$ satisfies the following stochastic evolution equation :
\begin{equation}\label{phi equation}
  d\phi_t(x)=\no(x,t)\left(-H_f(x,t)dt+\circ dW(t)\right),\quad x\in M.
\end{equation}
It should be emphasized that the choice of Stratonovich differential rather that It\^{o} differential is necessary to preserve geometric meaning of the equation that the equation is invariant under diffeomorphisms.
We restrict ourselves to the case of graphs over flat torus $\T^n$(i.e. the unit cube with periodic boundary condition), our random hypersurfaces are given by 
\[
  M(t,\omega)=\text{graph of }u(\cdot,t,\omega)=\left\{(x,u(x,t,\omega))\in\R^{n+1}:x\in\T^n,t\geq 0,\omega\in\Omega\right\}.
\]
In this setting, the outer normal $\no$ and the mean curvature $H$ are given by
\[
  \no=\frac{(-\n u,1)^{T}}{\sqrt{1+|\n u|^2}}, \quad H=\mathrm{div} \no=-\mathrm{div}\left(\frac{\n u}{\sqrt{1+|\n u|^2}}\right).
\]
The $f$-mean curvature is given by 
\begin{align*}
  H_f&=H-\langle \overline{\n} f,\no\rangle=-\mathrm{div}\left(\frac{\n u}{\sqrt{1+|\n u|^2}}\right)-\left\langle\overline{\n} f(x,u(x)),\frac{(-\n u,1)^{T}}{\sqrt{1+|\n u|^2}}\right\rangle\\
  &=-\mathrm{div}\left(\frac{\n u}{\sqrt{1+|\n u|^2}}\right)+\frac{\n f(x,u)\cdot \n u-\partial_{n+1}f(x,u)}{\sqrt{1+|\n u|^2}}.
\end{align*}
where $\overline{\n}f(x,x_{n+1})=(\n f,\partial_{n+1}f)^{\top}(x,x_{n+1})$, here that $\overline{\n}$ and $\n$ denote the derivative in $\R^{n+1}$ and $\R^n$ respectively,
for all $(x,x_{n+1})\in\R^{n+1}$. Suppose that $f$ satisfies
\[
  \n f(x,x_{n+1})=\n f(x),\quad \partial_{n+1}f(x,x_{n+1})=\xi x_{n+1},
\]
i.e. $f(x,x_{n+1})=f(x)+\frac{\xi}{2}x^2_{n+1}$, where $\xi $ is a constant. That is to say, $\n f$ depends only on the first $n$ variables $x$, and $\partial_{n+1}f$ depends only on the last variable $x_{n+1}$. It is easy to see that there exist many functions satisfying this condition, for example, $f(x)=c|x|^2$. Then, 
\[
  H_f=-\mathrm{div}\left(\frac{\n u}{\sqrt{1+|\n u|^2}}\right)+\frac{\n f\cdot \n u-\xi u}{\sqrt{1+|\n u|^2}}.
\]
then equation (\ref{phi equation}) reduces to the SPDE 
\begin{equation}
  du=\left(\sqrt{1+|\n u|^2}\mathrm{div}\left(\frac{\n u}{\sqrt{1+|\n u|^2}}\right)-\n f\cdot\n u+\xi u\right)dt+\sqrt{1+|\n u|^2}\circ dW.
\end{equation}
To simplify the notation, we define the divergence $\D_f V=\D V-\n f\cdot V$ for any vector field $V\in\R^n$ with respect to the weighted measure $d\mu=e^{-f}dx$, we denote $\Q, \V$ by 
\[
  \mathbf{Q}(p):=\sqrt{1+|p|^2},\quad \mathbf{v}(p):=\frac{p}{\sqrt{1+|p|^2}}, \quad p\in\R^n.
\]
Note that $\Q$ is the area element and $\V$ is the horizontal projection of the inner normal $-\no$ to the graph. With these notations, we have the simpler version SPDE 
\begin{equation}\label{SfMCF}
  du=\left(\Q\D_f\V+\xi u\right)dt+\Q\circ dW.  \tag{SFMCF}
\end{equation}

We will use the abstract theory of variational SPDEs to handle this equation. Since (\Ref{SfMCF}) lacks of coercivity,  an approximation with a viscous term $\epsilon\Delta_fu$ is made by
\begin{equation}\label{viscous equation}
  du=\left(\epsilon\Delta_fu+\Q\D_f\V+\xi u\right)+\Q\circ dW.
\end{equation}
where $\Delta_f u=\D_f \n u=\Delta u-\n f\cdot\n u$ is called the weighted Laplacian with respective to $d\mu=e^{-f(x)}dx,x\in \R^{n}$. The equation with a viscous term is coercive in some sense, we call it the viscous equation.

The classical theory \cite{Prévôt2007} for variational SPDEs dose not work for the viscous equation. We will adopt the compact embedding method developed in \cite{Dabrock2021}*{Appendix A}, which can be viewed as a generalization of the compact method developed by Pardoux \cite{Prévôt2007} and Viot \cite{Viot}.

This paper is organized as follows: In Sect. \ref{sec:2}, we state some ingredients that will be used in the following sections. In Sect. \ref{sec:3}, we list our main result. In sect. \ref{sec:4}, we prove the existence of martingale solution of the viscous equation \eqref{viscous equation}. In Sect. \ref{sec:5}, we prove a uniform estimate for the viscous equation \eqref{viscous equation}. In Sect. \ref{sec:6}, we give the proof of our main result. In Sect. \ref{sec:7}, we study the small perturbation limit of the stochastic $f$-mean curvature flow.

\section{Preliminary}\label{sec:2}
In this section we first introduce some basic notations and two useful formulas that will be used frequently in this paper, then we state the results in \cite{Dabrock2021} for the convenience of readers since we will use the results to prove our main theorem.

\subsection*{Periodic $L^p$ and Sobolev spaces}
Throughout this paper, the spaces we consider are $L^p$ or Sobolev spaces on $\T^n$ with the weighted measure $d\mu=e^{-f(x)}dx$ where $dx$ is the usual Lebesgue measure on $\T^n$. We note that the function $f$ in the weighted measure is a bit of difference with the $f$ in the $f$-mean curvature flow, the former is defined on $\R^n$ and the latter is defined on $\R^{n+1}$, the connection between them is $f(x,x_{n+1})=f(x)+\psi(x_{n+1})$, that is to say that the former is the horizontal projection of the restriction of the latter on the first $n$ variables $x\in \R^n$. For $k\geq 0, p\in [1,\infty]$ we denote $W^{k,p}(\T^n,d\mu)$ the space of periodic  Sobolev functions on the flat torus $\T^n$ which is the completion of the space of smooth periodic functions on $[0,1]^n$ with respect to the norm $\|\cdot\|_{W^{k,p}(\T^n,d\mu)}$ for $p<\infty$; and $W^{k,\infty}(\T^n,d\mu):=\{u\in W^{k,1}(\T^n,d\mu): |D^ju|\in L^{\infty}(\T^n,d\mu), 0\leq j\leq k\}$.

\subsection{Two useful formulas}
Now we prove two formulas that will be used frequently in our sequential proofs:
\begin{proposition}
  For $u,v\in C^2(\T^n)$, we have:
  \begin{align}
    \int_{\mathbb{T}^n}\nabla u\cdot\nabla vd\mu&=-\int_{\mathbb{T}^n}v\Delta_fud\mu,\\
    \int_{\mathbb{T}^n}|D^2u|^2d\mu&=\int_{\mathbb{T}^n}(\Delta_fu)^2d\mu-\int_{\mathbb{T}^n}\nabla u\cdot D^2f\nabla ud\mu.
  \end{align}
\end{proposition}
\begin{proof}
  For the first one: 
\begin{align*}
  \int_{\mathbb{T}^n}\nabla u\cdot\nabla vd\mu &= \sum_{i=1}^n\int_{\mathbb{T}^n}\partial_i u\partial_i ve^{-f}dx=-\int_{\mathbb{T}^n}\partial_i(\partial_i ue^{-f})vdx\\
  &=-\sum^n_{i=1}\int_{\mathbb{T}^n}(\partial_{ii}u-\partial_i f\partial_i u)ve^{-f}dx=-\int_{\mathbb{T}^n}(\Delta u-\nabla f\cdot\nabla u)vd\mu\\
  &=-\int_{\mathbb{T}^n}v\Delta_fud\mu.
\end{align*}
The second one:
\begin{align*}
  \int_{\mathbb{T}^n}|D^2u|^2d\mu &=\sum^n_{i=1}\int_{\mathbb{T}^n}\nabla\partial_i u\cdot\nabla\partial_i ud\mu=-\sum^n_{i=1}\int_{\mathbb{T}^n}\partial_i u\Delta_f\partial_iud\mu\\
  &=-\sum^n_{i=1}\int_{\mathbb{T}^n}\partial_iu(\partial_i\Delta_fu+\nabla\partial_if\cdot\nabla u)d\mu\\
  &=-\int_{\mathbb{T}^n}\nabla u\cdot\nabla\Delta_fud\mu-\int_{\mathbb{T}^n}\nabla u\cdot D^2f\nabla ud\mu\\
  &=\int_{\mathbb{T}^n}(\Delta_fu)^2d\mu-\int_{\mathbb{T}^n}\nabla u\cdot D^2f\nabla ud\mu.
\end{align*}
\end{proof}

\subsection{Variational SPDE under a compactness assumption}
The following results are taken from \cite{Dabrock2021} which will be used to prove our main result.
\begin{assumption}\label{ass1}
 Let $V$ and $H$ be separable Hilbert spaces with $V \subset H \simeq H^{\prime} \subset V^{\prime}$ and $V$ densely and compactly embedded in $H$. Furthermore, we will consider another separable Hilbert space $U$, which will be the space where a Wiener process is defined. Besides, we will denote by $\left(g_l\right)_{l \in \mathbb{N}}$ an orthonormal basis of $U$.
If not otherwise specified then a cylindrical Wiener process $W$ on $U$ with respect to a filtration $\left(\mathcal{F}_t\right)_t$ will always be assumed to have the representation $W=\sum_{l \in \mathbb{N}} g_l \beta_l$ with $\left(\beta_l\right)_{l \in \mathbb{N}}$ mutually independent real-valued $\left(\mathcal{F}_t\right)$-Brownian motions.

With $L_1(H)$ we will denote the space of all nuclear operators $T: H \rightarrow H$ with the norm
$$
\|T\|_{L_1(H)}:=\inf \left\{\sum_{k=1}^{\infty}\left\|a_k\right\|_H\left\|\varphi_k\right\|_H \mid\left(a_k\right)_k \subset H,\left(\varphi_k\right)_k \subset H^{\prime}, T=\sum_{k=1}^{\infty} a_k \varphi_k\right\} .
$$
It is well known that $\left(L_1(H)\right)^*=L(H)$ and that the weak-* topology on $L(H)$ coincides on norm bounded subsets with the weak operator topology on $L(H)$, which is the weakest topology such that for all $x, y \in H$ the map $L(H) \rightarrow \mathbb{R}, T \mapsto\langle T x, y\rangle_H$ is continuous.
\end{assumption}

N. Dabrock et al \cite{Dabrock2021}. consider the following SPDE 
\begin{equation}\label{spde}
  \begin{aligned}
    \mathrm{d} u & =A(u) \mathrm{d} t+B(u) \mathrm{d} W, \\
    u(0) & =u_0 .
  \end{aligned}
\end{equation}
They proved a generalized It\^{o} formula for variational SPDEs. 
\begin{proposition}[Itô formula and continuity]\label{Ito formula}
  Assume that $T>0,\left(\Omega, \mathcal{F},\left(\mathcal{F}_t\right)_{t \in[0, T]}, \mathbb{P}\right)$ is a stochastic basis with a normal filtration and $W$ a cylindrical Wiener process on $U$. Furthermore, let $u_0 \in L^2(\Omega ; H)$ be $\mathcal{F}_0$-measurable and $u, v, B$ be predictable processes with values in $V, V^{\prime}$ and $L_2(U ; H)$, respectively, such that
$$
u \in L^2\left(\Omega ; L^2(0, T ; V)\right), v \in L^2\left(\Omega ; L^2\left(0, T ; V^{\prime}\right)\right), B \in L^2\left(\Omega ; L^2\left(0, T ; L_2(U ; H)\right)\right) \text {, }
$$
and
$$
u(t)-u_0=\int_0^t v(s) d s+\int_0^t B(s) d W_s \text { in } V^{\prime} \mathbb{P} \text {-a.s. } \forall t \in[0, T].
$$
Then $u$ has a version with continuous paths in $H$ and for this version it holds that $u \in L^2(\Omega ; C([0, T] ; H))$ with
$$
\|u(t)\|_H^2-\left\|u_0\right\|_H^2=\int_0^t 2\langle v(s), u(s)\rangle_{V^{\prime}, V}+\|B(s)\|_{L_2(U ; H)}^2 d s+2 \int_0^t\langle u(s), B(s) d W_s\rangle_H \forall t \in[0, T].
$$

Furthermore, if $F \in C^1(H)$ and the second Gâteaux derivative $\mathrm{D}^2 F: H \rightarrow L(H)$ exists with
\begin{itemize}
  \item $F, D F$ and $\mathrm{D}^2 F$ bounded on bounded subsets of $H$,
  \item $\mathrm{D}^2 F: H \rightarrow L(H)$ continuous from the strong topology on $H$ to the weak-* topology on $L(H)=\left(L_1(H)\right)^*$ and
  \item $\left.(D F)\right|_V: V \rightarrow V$ continuous from the strong topology on $V$ to the weak topology on $V$ and growing only linearly
$$
\|D F(x)\|_V \leq C\left(1+\|x\|_V\right) \forall x \in V,
$$
\end{itemize}
then $\mathbb{P}$-a.s. for all $t \in[0, T]$
$$
\begin{aligned}
F(u(t))-F\left(u_0\right)= & \int_0^t\langle v(s), D F(u(s))\rangle_{V^{\prime}, V}+\frac{1}{2} \operatorname{tr}\left[\mathrm{D}^2 F(u(s)) B(s)(B(s))^*\right] d s \\
& +\int_0^t\langle D F(u(s)), B(s) d W_s\rangle_H .
\end{aligned}
$$
\end{proposition}

In addition to Assumption \ref{ass1}, N. Dabrock et al. \cite{Dabrock2021} make the following assumptions.
\begin{assumption}\label{ass2}
  Let $A: V \rightarrow V^{\prime}$ and $B: V \rightarrow L_2(U ; H)$. We will write $B^*: V \rightarrow$ $L_2(H ; U)$ for the adjoint operator $B^*(u):=(B(u))^*$. Assume:
  \begin{itemize}
    \item \textbf{Coercivity: } There exist constants $\alpha, C>0$ such that
        \begin{equation*}
          2\langle A(u), u\rangle_{V^{\prime}, V}+\|B(u)\|_{L_2(U ; H)}^2 \leq-\alpha\|u\|_V^2+C\left(1+\|u\|_H^2\right) \forall u \in V .
        \end{equation*} 
    \item \textbf{Growth bounds: } There is a constant $C>0$ and $\delta \in(0,2]$ such that \begin{align*}
          \|A(u)\|_{V^{\prime}}^2 & \leq C\left(1+\|u\|_V^2\right) \quad\forall u \in V, \\
          \|B(u)\|_{L_2(U ; H)}^2 & \leq C\left(1+\|u\|_V^2\right) \quad \forall u \in V, \\
          \|B(u)\|_{L\left(U ; V^{\prime}\right)}^2 & \leq C\left(1+\|u\|_V^{2-\delta}+\|u\|_H^2\right) .
        \end{align*}
    \item \textbf{Continuity: } $A: V \rightarrow V^{\prime}$ is weak-weak- $\ast$  sequentially continuous, that means\[
      u_k \rightarrow u \text { in } V \Rightarrow A\left(u_k\right) \stackrel{*}{\rightharpoonup} A(u) \text { in } V^{\prime},
    \]
    and $B^*: V \rightarrow L_2(H ; U)$ is sequentially continuous from the weak topology on $V$ to the strong operator topology on $L(H ; U)$, that means\[
      u_k \rightarrow u \text { in } V \Rightarrow B^*\left(u_k\right) h \rightarrow B^*(u) h \,\text { in } U \,\forall h  \in H.
    \]
  \end{itemize}
\end{assumption}
Under the Assumption \ref{ass1} and Assumption \ref{ass2}, N. Dabrock et al. could get the existence of martingale solution of (\ref{spde}). They first apply classical theory for finite dimensional SDEs to prove the existence weak martingale solutions of finite dimensional Galerkin approximations of (\ref{spde}). Then the It\^{o} formula under the Assumption \ref{ass2} is used to deduce uniform bound of $\|u(t)\|^2_H$. After that, they use the Jakubowski-Skorokhod representation theorem for tight sequences in non-metric spaces to prove the approximations converge on some probability space. Finally they show that this limit is exactly the martingale solution of (\ref{spde}).
\begin{theorem}\label{existence}
   Let $q>2$ and $\Lambda$ be a Borel probability measure on $H$ with finite $q$-th moment
$$
\int_H\|z\|_H^q d \Lambda(z)<\infty.
$$
Then there is a martingale solution of (\ref{spde}) with initial data $\Lambda$. That means, that there is a stochastic basis $\left(\Omega, \mathcal{F},\left(\mathcal{F}_t\right)_{t \in[0, \infty)}, \mathbb{P}\right)$ with a normal filtration, a cylindrical $\left(\mathcal{F}_t\right)$-Wiener process $W$ on $U$ and a predictable $u$ with $u \in L^2\left(\Omega ; L^2(0, T ; V)\right) \cap$ $L^2(\Omega ; C([0, T] ; H))$ for all $T>0$ and
$$
\begin{gathered}
\langle u(t), v\rangle_H-\langle u(0), v\rangle_H=\int_0^t\langle A(u(s)), v\rangle_{V^{\prime}, V} d s+\int_0^t\langle B(u(s)) d W_s, v\rangle_H \\
=\int_0^t\langle A(u(s)), v\rangle_{V^{\prime}, V} d s+\sum_{l \in \mathbb{N}} \int_0^t\left\langle B(u(s)) g_l, v\right\rangle_H d \beta_l(s).
\end{gathered}
$$
$\mathbb{P}$-a.s. for all $t \in[0, \infty)$ and $v \in V$, and $\mathbb{P} \circ u(0)^{-1}=\Lambda$.
\end{theorem}

\section{Results}\label{sec:3}
In the section we state our main results of this paper. First we formulate the concept of a martingale solution of a SPDE.

\begin{definition}\label{definition}
  \begin{enumerate}
    \item [(i)] Let $I=[0,\infty)$ or $I=[0,T]$ for some $T>0$. $(\Omega,\mathcal{F},(\mathcal{F}_t)_{t\in I},\p)$ be a stochastic basis with a normal filtration together with a real-valued $(\mathcal{F}_t)$-Wiener process $W$ and $u_0\in L^2(\Omega; H^1(\T^n,d\mu))$ be $\mathcal{F}_0$-measurable. A predictable $H^2(\T^n,d\mu)$-valued process $u$ with $u\in L^2(\Omega; L^2(0,T; H^2(\T^n,d\mu)))$ for all $t\in I$ is a strong solution of (\ref{SfMCF}) with initial data $u_0$, if \[
      u(t)-u_0=\int^t_0\Q\D_f\V +\xi uds+\int^t_0\Q\circ dW_s\quad \p\text{-a.s. in}\, L^2(\T^n,d\mu)\, \forall t\in I.
    \]
    \item [(ii)] Let $\Lambda$ be a probability measure on $H^1(\T^n,d\mu)$ with bounded seconds moments, i.e. \[\int_{H^1(\T^n,d\mu)}\|z\|^2_{H^1}d\Lambda(z)<\infty.\] A martingale solution of (\ref{SfMCF}) with initial date $u_0$ is given by $(\Omega,\mathcal{F},(\mathcal{F}_t)_{t\in I},\p)$ together with $W, u_0$ and $u$ such that (i) is satisfied and $\p\circ u^{-1}_0=\Lambda$.
  \end{enumerate}
\end{definition}
We can also define the strong and martingale solutions for the viscous equation (\ref{viscous equation}). In the following we will just write $u$ is a martingale solution to simply the notation of martingale solution $(u,(\Omega,\mathcal{F},(\mathcal{F}_t)_{t\in I},\p),W,u_0)$ if there is no danger of confusion. And the term `solution' we use refers to martingale solution.

Formally for a solution of (\ref{SfMCF}), we can use chain rule for Stratonovich differential to deduce that 
\[
  d\Q=\ast \ast dt+\V\cdot\n \left(\Q\right)\circ dW.
\]
Hence, the It\^{o}-Stratonovich correlation for the integral in Definition \ref{definition} is given by 
\[
  \frac{1}{2}\V\cdot \n \left(\Q\right)=\frac{1}{2}\V\cdot D^2u\V,
\]
the Stractonovich integral in Definition \ref{definition} becomes
\[
  \int^t_0\Q\circ dW_s=\int^t_0\Q dW_s+\frac{1}{2}\int^t_0\V\cdot D^2u\V ds.
\] 
Thus, the equation in Definition \ref{definition} becomes
\begin{equation}
  \begin{aligned}
    u(t)-u_0&=\int^t_0\left[\Q\D_f\left(\V\right)+\xi u+\frac{1}{2}\V\cdot D^2u\V\right]dt+\int^t_0\Q dW_s\\
    &=\int^t_0\left(\Delta_fu-\frac{1}{2}\V \cdot D^2u\V+\xi u\right)ds+\int^t_0\Q dW_s.
  \end{aligned}  
\end{equation}
Now we state our main result of the paper.
\begin{theorem}[Existence of martingale solutions]\label{main theorem}
  Let $\Lambda$ be a Borel probability measure on $H^1(\T^n,d\mu)$ with bounded second moments and additionally
  \[
    \mathrm{supp} \Lambda\subset \{z\in H^1(\T^n,d\mu): \|\nabla z\|_{L^\infty}\leq L\},
  \] for some constant $L>0$.

  Then for $I=[0,\infty)$, under the pinching condition: there is a constant $\delta>0$ such that \[D^2f\geq \delta\xi\mathrm{Id}> 0,\] 
  there is a martingale solution of (\ref{SfMCF}) with initial date $\Lambda$. For all such solutions it holds that $D^2u\in L^2(\Omega; L^2(0,\infty;L^2(\T^n,d\mu)))$ and 
  \[
    \|\n u\|_{L^{\infty}(0,\infty;L^{\infty})}\leq L\quad \p-a.s..
  \]
\end{theorem}

\section{Existence of martingale solutions of viscous equation}\label{sec:4}
In this section, we will use the theory listed in Section 2.2 or in \cite{Dabrock2021}*{Appendix A} to prove the existence of martingale solution of the viscous equation (\ref{viscous equation}) of the stochastic $f$-mean curvature flow. We apply the framework listed in Section 2.2 to $\n u$ rather than directly to $u$. Then we can construct $u$ by $\n u$ and it is the solution of the viscous equation (\ref{viscous equation}).

\begin{theorem}\label{existence of viscous martingale solution}
  Let $\epsilon>0, q>2$ and $\Lambda$ be a Borel probability measure on $H^1(\T^n,d\mu)$ with 
  \[
    \int_{H^1(\T^n,d\mu)}\|z\|^2_{H^1}d\Lambda(z)<\infty,
  \]
  and
  \[
    \int_{H^1(\T^n,d\mu)}\|\n z\|_{L^2,d\mu}^q<\infty.
  \]
  Under condition: \[D^2f\geq 0,\quad \xi\geq 0.\] 
  there is a martingale solution $u$ of the viscous equation(\ref{viscous equation}) for $I=[0,\infty)$ with initial data $\Lambda$.
\end{theorem}
\begin{proof}
We first consider the SPDE that $\n u$ satisfies for the viscous equation:
\begin{equation}\label{gradient}
  d\n u=\n \left(\epsilon\Delta_f+\Q\D_f\V+\frac{1}{2}\V\cdot D^2u\V+\xi u\right)+\n \Q dW.
\end{equation}
We will work with the spaces:
$$
\begin{aligned}
  &V:=\left\{\nabla u | u\in H^2(\mathbb{T}^n,d\mu)\right\} \text{with} \|\nabla u\|_V:=\|\nabla u\|_{H^1(\mathbb{T}^n;\mathbb{R}^n,d\mu)};\\
  &H:=\left\{\nabla u | u\in H^1(\mathbb{T}^n,d\mu)\right\} \text{with} \|\nabla u\|_H:=\|\nabla u\|_{L^2(\mathbb{T}^n;\mathbb{R}^n,d\mu)};\\
  &U:=\mathbb{R}.
\end{aligned}
$$
We have $V\subset H$ densely and compactly. Furthermore, we can identify $L_2(U;H)=H.$
We define the operators:
\[
  \begin{aligned}
    A_{\epsilon}&: V\longrightarrow V^{\prime}\\
    \langle A_{\epsilon}(\nabla u),\nabla w\rangle_{V^{\prime},V}
    &:=\int_{\mathbb{T}^n}\nabla \left(\epsilon\Delta_f+\Q\D_f\V+\frac{1}{2}\V\cdot D^2u\V+\xi u\right)\cdot\nabla wd\mu\\
    &=-\int_{\mathbb{T}^n}\left((1+\epsilon)\Delta_f u-\frac{1}{2}\V D^2u\V\right)\Delta_f wd\mu+\int_{\T^n}\xi \n u\cdot\n wd\mu.
  \end{aligned}
\]
and
\[
  \begin{aligned}
    &B:V\longrightarrow H\\
    &B(\nabla u):=\nabla\Q=D^2u \V.
  \end{aligned}
\]
Since $D^2f=(\partial_{ij}f)$ is convex and continuous
on $\mathbb{T}^n$, there exists a uniform constant $M$ such that $\sup_{\mathbb{T}^n}|\partial_{ij}f|\leq M$ for all $i,j=1,2,\dots,n$. Then we have 
\[
  0\leq \nabla u\cdot D^2f\nabla u=\sum_{ij}\partial_iu\partial_{ij}f\partial_ju\leq M\left(\sum_i\partial_i u\right)^2\leq nM|\nabla u|^2.
\]

We now verify the assumption are fulfilled.

 $\bullet$ \textbf{Coercivity}
       \[
        \begin{aligned}
          2&\langle A_{\epsilon}(\nabla u),\nabla u\rangle_{V^{\prime},V}+\|B(\nabla u)\|^2_H\\
          &=\int_{\mathbb{T}^n}-2\left((1+\epsilon)\Delta_f u-\frac{1}{2}\V D^2u\V\right)\Delta_f u+2\xi|\n u|^2d\mu+\int_{\T^n}\left|D^2u\V\right|^2d\mu\\
          &=-2(1+\epsilon)\int_{\mathbb{T}^n}(\Delta_f u)^2d\mu+\int_{\mathbb{T}^n}\V\cdot D^2u\V \Delta_f ud\mu+\int_{\mathbb{T}^n}2\xi|\n u|^2+\left|D^2u\V\right|^2d\mu\\
          &\leq -2(1+\epsilon)\int_{\mathbb{T}^n}(\Delta_f u)^2d\mu+\frac{1}{2}\int_{\mathbb{T}^n}|D^2u|^2+(\Delta_f u)^2d\mu+\int_{\mathbb{T}^n}2\xi|\n u|^2+|D^2u|^2d\mu\\
          &=-2\epsilon\int_{\mathbb{T}^n}(\Delta_f u)^2d\mu-\frac{3}{2}\int_{\mathbb{T}^n}\nabla u\cdot D^2f\nabla ud\mu+2\xi\int_{\T^n}|\n u|^2d\mu\\
          &=-2\epsilon\int_{\mathbb{T}^n}|D^2u|^2d\mu-\left(2\epsilon+\frac{3}{2}\right)\int_{\mathbb{T}^n}\nabla u\cdot D^2f\nabla ud\mu+2\xi\int_{\T^n}|\n u|^2d\mu\\
          &\leq-2\epsilon\|\nabla u\|^2_{H^1(\mathbb{T}^n;\mathbb{R}^n)}+2\left(\epsilon+\xi\right)\|\nabla u\|^2_{L^2({\mathbb{T}^n};\mathbb{R}^n)}.
         \end{aligned}
       \]
       Coercivity is satisfied for  $\epsilon>0$.

   $\bullet$ \textbf{Growth bounds}
   \begin{align*}
    \|A_{\epsilon }&(\nabla u)\|^2_{V^{\prime}}\\
    & =\left(\sup_{\nabla g\in V}\frac{\int_{\mathbb{T}^n}\nabla\left((1+\epsilon)\Delta_f u-\frac{1}{2}\V\cdot D^2u\V +\xi u\right)\cdot\nabla gd\mu}{\|\nabla g\|_{H^1(\mathbb{T}^n;\mathbb{R}^n,d\mu)}}\right)^2\\
    & =\sup_{\nabla g\in V}\frac{\left(-\int_{\mathbb{T}^n}\left((1+\epsilon)\Delta_f u-\frac{1}{2}\V\cdot D^2u\V\right)\Delta_f gd\mu+\xi\int_{\T^n}\n u\cdot\n gd\mu\right)}{\|\nabla g\|^2_{H^1(\mathbb{T}^n;\mathbb{R}^n,d\mu)}}^2\\
    & \leq\sup_{\nabla g\in V}\frac{2\int_{\mathbb{T}^n}(\Delta_f g)^2d\mu}{\|\nabla g\|^2_{H^1(\mathbb{T}^n;\mathbb{R}^n,d\mu)}}\cdot\int_{\mathbb{T}^n}\left((1+\epsilon)\Delta_f u-\frac{1}{2}\V\cdot D^2u\V\right)^2d\mu\\
    &\qquad\quad +\sup_{\n g\in V}\frac{2\xi^2\int_{\T^n}|\n u|^2d\mu\int_{\T^n}|\n g|^2d\mu}{\|\n g\|^2_{H^1_{\T^n;\R^n,d\mu}}}\\
    & \leq\sup_{\nabla g\in V}\frac{2\int_{\mathbb{T}^n}|D^2g|^2+\nabla g\cdot D^2f\nabla gd\mu}{\int_{\mathbb{T}^n}|D^2 g|^2+|\nabla g|^2d\mu}\cdot \int_{\mathbb{T}^n}2(1+\epsilon)^2(\Delta_f u)^2+\frac{1}{2}|D^2u|^2d\mu\\
    &\qquad\quad +\sup_{\n g\in V}\frac{2\xi^2\int_{T^n}|\n u|^2d\mu\int_{T^n}|\n g|^2d\mu}{\int_{\T^n}|\n g|^2+|D^2g|^2d\mu}\\
    & \leq\sup_{\nabla g\in V}\frac{2\int_{\mathbb{T}^n}|D^2 g|^2+nM|\nabla g|^2d\mu}{\int_{\mathbb{T}^n}|D^2 g|^2+|\nabla g|^2d\mu}\\
    &\quad\quad\quad\cdot\int_{\mathbb{T}^n}\left[2(1+\epsilon)^2+\frac{1}{2}\right]|D^2u|^2+2(1+\epsilon)^2\nabla u\cdot D^2f\nabla ud\mu+2\xi^2\int_{\T^n}|\n u|^2d\mu\\
    & \leq \int_{\mathbb{T}^n}2(1+nM)\left[2(1+\epsilon)^2+\frac{1}{2}\right]|D^2u|^2+\left(4nM(1+nM)(1+\epsilon)^2+2\xi^2\right)|\nabla u|^2d\mu\\
    & \leq \left(2(1+nM)\left[2(1+nM)(1+\epsilon)^2+\frac{1}{2}\right]+2\xi^2\right)\|\nabla u\|^2_{H^1(\mathbb{T}^n;\mathbb{R}^n,d\mu)}.\\
    \|B(\nabla u)\|^2_H
    &=\|D^2u\V\|^2_H\leq \|D^2u\|^2_{L^2}\leq\|\nabla u\|^2_{H^1}=\|\nabla u\|^2_V.\\
    \|B(\nabla u\|)^2_{V^{\prime}}
    &=\left(\sup_{\nabla g\in V}\frac{\int_{\mathbb{T}^n}\nabla \Q\cdot\nabla gd\mu}{\|\nabla g\|_{H^1}}\right)^2\\
    &\leq 2(1+nM)\int_{\mathbb{T}^n}\Q^2d\mu.\\
    &=2(1+nM)\int_{\T^n}1+|\n u|^2d\mu\\
    &\leq 2(1+nM)\left(1+\|\nabla u\|^2_{L^2(\mathbb{T}^n;\mathbb{R}^n,d\mu)}\right).\notag
   \end{align*}

 $\bullet$ \textbf{Continuity}:
       When $\nabla u_k\rightharpoonup \nabla u$ in $V$, then $\nabla u_k\to \nabla u$ in $H$ and $\frac{\nabla u_k}{|\nabla u_k|}\rightharpoonup\frac{\nabla u}{|\nabla u|}$ in $H$. Therefore,
       \begin{align*}
        \mathbf{v}\left(\nabla u_k\right) \cdot \mathrm{D}^2 u_k \mathbf{v}\left(\nabla u_k\right) & =\mathbf{v}\left(\nabla u_k\right) \otimes \mathbf{v}\left(\nabla u_k\right): \mathrm{D}^2 u_k \\
        & \rightarrow \mathbf{v}(\nabla u) \otimes \mathbf{v}(\nabla u): \mathrm{D}^2 u \\
        & =\mathbf{v}(\nabla u) \cdot \mathrm{D}^2 u \mathbf{v}(\nabla u) \text { in } L^1\left(\T^n,d\mu\right).
      \end{align*}
      since $|\mathbf{v}(\n u_k)|\leq 1$, we have 
      \[
        \mathbf{v}\left(\nabla u_k\right) \cdot \mathrm{D}^2 u_k \mathbf{v}\left(\nabla u_k\right) \rightarrow \mathbf{v}(\nabla u) \cdot \mathrm{D}^2 u \mathbf{v}(\nabla u) \text { in } L^2\left(\mathbb{T}^n,d\mu\right).
      \]
      The other terms in the definition of $A(u_k)$ are linear in $u_k$, hence we have 
      \[
        A(u_k)\stackrel{\ast}{\rightharpoonup } A(u) \text{ in } V^{\prime}.
      \]
      Similarly, we have \[B(\nabla u_k)\rightharpoonup B(\nabla u) \text{ in } H.\]
Thus, Assumption \ref{ass2} is fulfilled, by Theorem \ref{existence} we get a martingale solution of (\ref{gradient}).

Next, we show this equation is fulfilled in $H^{-1}(\mathbb{T}^n;\mathbb{R}^n,d\mu),$ hence the solution is weak in PDE sense.
For arbitrary $\phi\in H^1(\mathbb{T}^n;\mathbb{R}^n,d\mu)$, with the help of the Helmholtz-Hodge decomposition theorem for weighted measure \cite{XD-Li}*{Theorem 1.9}, we can write $\phi=\nabla\psi+\omega$ with $\psi\in H^2(\mathbb{T}^n;\mathbb{R}^n,d\mu)$ and $\operatorname{div}_f \omega:=\operatorname{div}\omega -\nabla f\cdot\omega=0$, then we have
$$
\begin{aligned}
\int_{\mathbb{T}^n}(\nabla u(t)-\nabla u_0)\cdot\phi d\mu
&=\int^t_0\langle A_{\epsilon}(\nabla u(s)),\phi\rangle_{H^{-1},H^1}ds+\int^t_0\int_{\mathbb{T}^n}B(\nabla u(s))\cdot\phi d\mu dW_s.
\end{aligned}
$$
Therefore, $\nabla u$ is fulfilled in $H^{-1}(\mathbb{T}^n;\mathbb{R}^n,d\mu).$

Now, we define for $t\in [0,\infty)$,
\[
  F(t):=u_0+\int^t_0\left((1+\epsilon)\Delta_f u(s)-\frac{1}{2}\V D^2u\V(s)\right)ds+\int^t_0\Q dW_s.
\]
With the assumption that $u_0\in L^2(\Omega;L^2(\mathbb{T}^2,d\mu))$, and for $T\in [0,\infty)$, it is easy to check that 
$$
\begin{aligned}
  t & \longmapsto  \int^t_0\left((1+\epsilon)\Delta_f u(s)-\frac{1}{2}\V D^2u(s)\V \right)dt \in L^2(\Omega;L^2(0,T;L^2(\mathbb{T}^n,d\mu))),\\
  t & \longmapsto  \int^t_0\Q dW_s\in L^2(\Omega;L^2(0,T;L^2(\mathbb{T}^n,d\mu))).
\end{aligned}
$$
Thus, $F(t)\in L^2(\Omega;L^2(0,T;L^2(\mathbb{T}^n,d\mu)))$.

Since $\n u(t)$ is a solution of (\ref{viscous equation}), by definition we have 
\begin{align*}
  \n u(t)-&\xi\int^t_0\n u(s)ds\\
  &=\n u_0+\int^t_0\n\left((1+\epsilon)\Delta_f u(s)-\frac{1}{2}\V D^2u\V(s)\right)ds+\int^t_0D^2u\V(s) dW_s\\
  &=\n F(t),
\end{align*}
thus, we get 
\[
  \n u(t)=\n F(t)+\xi\int^t_0e^{\xi(t-s)}\n F(s)ds.
\]

Next, we define $t\in [0,\infty)$,
\[
  \tilde{u}(t):=F(t)+\xi\int^t_0e^{\xi(t-s)}F(s)ds.
\]
Clearly, \[\n \tilde{u}(t)=\n u(t) \quad  \forall t\in[0,\infty)\, \p-a.s.\] It is easy to see that $\tilde{u}(0)=u_0$ and $\tilde{u}(t)\in L^2(\Omega;L^2(0,T;L^2(\mathbb{T}^n,d\mu)))$, what is more,
\[
  \int^t_0\tilde{u}(s)=\int^t_0e^{\xi(t-s)}F(s)ds.
\]
Thus, 
\begin{align*}
  \tilde{u}(t)&=F(t)+\xi\int^t_0e^{\xi(t-s)}F(s)ds\\
  &=u_0+\int^t_0\left((1+\epsilon)\Delta_f u(s)-\frac{1}{2}\V D^2u\V(s)\right)ds+\int^t_0\Q dW_s+\xi\int^t_0\tilde{u}(s)ds\\
  &=u_0+\int^t_0\left((1+\epsilon)\Delta_f \tilde{u}(s)-\frac{1}{2}\mathbf{v}(\n\tilde{u})D^2\tilde{u}\mathbf{v}(\n\tilde{u})(s)+\xi \tilde{u}(s)\right)+\int^t_0\mathbf{Q}(\n \tilde{u})(s)dW_s\\
  &=\left(\epsilon\Delta_f\tilde{u}+\mathbf{Q}(\n\tilde{u})\D_f\mathbf{v}(\n\tilde{u})+\xi \tilde{u}\right)+\mathbf{Q}(\n\tilde{u})\circ dW.
\end{align*}
Which means that $\tilde{u}$ is a martingale solution of the viscous equation (\ref{viscous equation}).

\end{proof}

\section{A priori energy estimates }\label{sec:5}
This section, we get a priori energy estimates for solutions of viscous equation (\ref{viscous equation}) which are uniformly in $\epsilon>0$ and also holds for true for solutions of the stochastic $f$-mean curvature flow (\ref{SfMCF}). This estimate is similar with \cite{Dabrock2021}*{Proposition 5.1} and can be viewed as its extension.
\begin{proposition}\label{uniform estimate of energy}
  Let $\epsilon\geq 0$ and $u$ be a solution of (\ref{viscous equation}). Under the pinching condition: there is a constant $\delta>0$ such that \[D^2f\geq \delta\xi\mathrm{Id}\geq 0,\] the energy $\int_{\mathbb{T}^n}|\nabla u|^2d\mu$ is a surpermartingale. 

  More explicitly, 
  \[
    \begin{aligned}
      \mathbb{E}\|\nabla u\|^2_{L^2(\mathbb{T}^2,d\mu)}&+2\epsilon\mathbb{E}\int^t_0\int_{\mathbb{T}^n}|D^2u(s)|^2d\mu ds+\frac{1}{2}\mathbb{E}\int^t_0\int_{\mathbb{T}^n}\mathbf{Q}(\n u(s))^2|\operatorname{div}_f v(s)|^2d\mu ds\\ 
      &+\E\int^t_0\int_{\T^n}\n u\cdot\left(\left(1+\frac{3\alpha}{2}\right)D^2f-2\alpha\xi\mathrm{Id}\right)\n ud\mu ds\\ 
      &+\mathbb{E}\int^t_0\int_{\mathbb{T}^n}\left(\frac{3}{2}|D^2u|^2-|D^2uv|^2-\frac{1}{2}|v\cdot D^2u v|^2\right)(s)d\mu ds\\ 
      \leq &\mathbb{E}\|\nabla u_0\|_{L^2(\mathbb{T}^n,d\mu)}.
    \end{aligned}
  \]
  where $\alpha$ satisfies $0<\dfrac{4\alpha}{2+3\alpha}\leq \delta$.
  We also have for $q\in [1,2)$ with a universal constant $C>0$, that 
  \[
    \mathbb{E}\sup_{t\in I}\|\nabla u\|^{2q}_{L^2(\T^n,d\mu)}\leq\left(2+\frac{2C^2}{2q-q^2}\right)\mathbb{E}\|\nabla u_0\|^{2q}_{L^2(\T^n,d\mu)}.
  \]
  If in addition $\p\text{-esssup}\|\n u\|_{L^{\infty}(I;L^{\infty}(\T^n))}=L<\infty$, then we have 
  \begin{align*}
    \mathbb{E}\|\nabla u(t)\|_{L^2\left(\mathbb{T}^n,d\mu\right)}^2 & +\frac{3+4 L^2}{2\left(1+L^2\right)^2} \mathbb{E} \int_0^t \int_{\mathbb{T}^n}\left|D^2 u(s)\right|^2 \mathrm{~d} s \leq \mathbb{E}\left\|\nabla u_0\right\|_{L^2\left(\mathbb{T}^n,d\mu\right)}^2 \quad \forall t \in I .
    \end{align*}
\end{proposition}

Actually, the uniform bounded assumption in Proposition \ref{uniform estimate of energy} of $\n u$ is fulfilled if the initial data $\n u_0$ is Lipschitz continuous with the same Lipschitz constant $L$.

\begin{proposition}[Maximum principle for the gradient of solutions]\label{Maximum principle}
  Let $\epsilon\geq 0$ and $u$ be a solution of (\ref{viscous equation}). 
  Under the same pinching condition as in Proposition \ref{uniform estimate of energy}, if $\mathbb{P}$-$\operatorname{esssup}\|\nabla u_0\|_{L^{\infty}(\mathbb{T}^n)}<\infty$, then 
  $\nabla u \in L^{\infty}(I;L^{\infty}(\mathbb{T}^n))$ a.s. with
  \[
    \|\nabla u\|_{L^{\infty}(I;L^{\infty}(\mathbb{T}^n))}\leq \mathbb{P}-\operatorname{esssup}\|\nabla u_0\|_{L^{\infty}(\mathbb{T}^n)}.
  \]
\end{proposition}
\begin{remark}
  For MCF, the same conclusion holds under the initial condition $\|\n u_0\|_{L^{\infty}}<L$. See \cite{Huisken1989}*{Corollary 3.2}.
\end{remark}

Both Proposition \ref{uniform estimate of energy} and Proposition \ref{Maximum principle} are proved based on It\^{o} formulas for integral of some special functions of the gradient of solutions.
\begin{lemma}\label{lemma of Ito formula of I}
  Let $\epsilon\geq 0$ and $u$ be the solution of (\ref{viscous equation}). For $g\in C^2(\mathbb{R}^n)$ with bounded 
  second order derivative and
  \[
    \mathcal{I}(t):=\int_{\mathbb{T}^n}g(\nabla u(t))d\mu, \quad t\in I,
  \]
  we have 
  \begin{equation}\label{Ito formula of I}
    \begin{aligned}
      d\mathcal{I}=
      &\int_{\mathbb{T}^n}-\epsilon D^2g(\nabla u)D^2u : D^2ud\mu+\int_{\mathbb{T}^n}-\left(\epsilon+\frac{1}{2}\right)\nabla g(\nabla u)\cdot D^2f\nabla ud\mu\\
      &+\int_{\mathbb{T}^n}-\frac{1}{2}g(\nabla u)|\operatorname{div}_f\mathbf{v}(\nabla u)|^2d\mu+\int_{\mathbb{T}^n}\frac{1}{2}\mathbf{v}(\nabla u)\cdot D^2f\left(g(\nabla u)\mathbf{v}(\nabla u)-\Q\nabla g(\nabla u)\right)d\mu\\
      &+\int_{\mathbb{T}^n}D^2u(\mathrm{Id}-\mathbf{v}(\nabla u)\otimes\mathbf{v}(\nabla u)):\left(\frac{g(\nabla u)}{2\Q ^2}(\mathrm{Id}-\mathbf{v}(\nabla u)\otimes \mathbf{v}(\nabla u))-D^2g(\nabla u)\right)D^2ud\mu\\
      &+\int_{\mathbb{T}^n}\frac{1}{2}(D\mathbf{v}(\nabla u))^{\top}:\left(\mathbf{v}(\nabla u)\otimes D^2u\nabla g(\nabla u)-\nabla g(\nabla u)\otimes D^2u\mathbf{v}(\nabla u)\right)d\mu dt\\
      &+\xi\int_{\T^n}\n g(\n u)\cdot\n ud\mu dt -\int_{\mathbb{T}^n}g(\nabla u)\operatorname{div}_f\mathbf{v}(\nabla u)d\mu dW.
    \end{aligned}  
  \end{equation}
\end{lemma}  
\begin{proof}
  To abbreviate the notation, we will write $Q=\Q, v=\mathbf{v}(\nabla u)$. Then we have $\nabla Q=D^2uv$. Applying
  the It\^{o} formula in Lemma \ref{Ito formula} infer 
  \begin{equation}
    \begin{aligned}
      d\mathcal{I}=
      &\int_{\mathbb{T}^n}\nabla g(\nabla u)\cdot\nabla \left(\epsilon\Delta_fu+Q\operatorname{div}_fv+\xi u+\frac{1}{2}v\cdot\nabla Q\right)d\mu dt\\
      &+\int_{\mathbb{T}^n}\frac{1}{2}\nabla Q\cdot D^2g(\nabla u)\nabla Qd\mu dt+\int_{\mathbb{T}^n}\nabla g(\nabla u)\cdot\nabla Qd\mu dW\\
      &=:\epsilon\lambda_{viscous}+\frac{1}{2}\lambda_{fmcf}+\frac{1}{2}\lambda_{pert} dt+\sigma dW.
    \end{aligned}
  \end{equation}
  with
  \[
    \begin{aligned}
    \lambda_{viscous}&=\int_{\mathbb{T}^n}\nabla g(\nabla u)\cdot\nabla\Delta_f u d\mu\\
    \lambda_{fmcf}&=\int_{\mathbb{T}^n}\nabla g(\nabla u)\cdot\nabla(Q\operatorname{div}_f v+\xi u)d\mu\\
    \lambda_{pert}&=\int_{\mathbb{T}^n}\nabla g(\nabla u)\cdot\nabla\Delta_f u +\nabla Q\cdot D^2g(\nabla u)\nabla Q  d\mu\\
    \sigma&=\int_{\mathbb{T}^n}\nabla g(\nabla u)\cdot\nabla Qd\mu.
  \end{aligned}
  \]
  The term $\lambda_{viscous}$ corresponds to the directional derivative of functional $u\mapsto \int_{\mathbb{T}^n}g(\nabla u)d\mu$ 
  into the direction of solutions of the weighted heat equation $\partial_t u=\Delta_f u$. The term $\lambda_{fmcf}$ corresponds to the directional derivative of functional $u\mapsto \int_{\mathbb{T}^n}g(\nabla u)d\mu$ 
  into the direction of solutions of unperturbed $f-$mean curvature flow of graphs. It is multiplied with $\frac{1}{2}$
  because the other half is used in $\lambda_{pert}$ to handle the additional terms coming from the perturbation. For $\lambda_{viscous}$ we calculate
  \begin{align*}
    \lambda_{viscous}
    &=\int_{\mathbb{T}^n}\nabla g(\nabla u)\cdot\nabla\Delta_f u d\mu=\int_{\mathbb{T}^n}\sum_{i}\partial_ig(\nabla u)\partial_i\left(\sum_j\partial_{jj}u-\partial_jf\partial_j u\right)d\mu\\
    &=\int_{\mathbb{T}^n}\sum_{ij}-\partial_j(\partial_i g(\nabla u))\partial_{ij}u+\partial_i g(\nabla u)\partial_{ij}u\partial_j fd\mu\\
    &\quad\quad-\int_{\mathbb{T}^n}\sum_{ij}\partial_i g(\nabla u)\partial_{ij}f\partial_j u+\partial_ig(\nabla u)\partial_jf\partial_{ij}ud\mu\\
    &=\int_{\mathbb{T}^n}-\sum_{ij}\nabla\partial_ig(\nabla u)\cdot\partial_j\nabla u\partial_{ij}u+\nabla g(\nabla u)D^2u\nabla fd\mu\\
    &\quad\quad-\int_{\mathbb{T}^n}\nabla g(\nabla u)D^2f\nabla u+\nabla g(\nabla u)D^2u\nabla fd\mu\\
    &=\int_{\mathbb{T}^n}-D^2g(\nabla u)D^2u:D^2ud\mu-\int_{\mathbb{T}^n}\nabla g(\nabla u)D^2f\nabla ud\mu.
  \end{align*}
  For $\lambda_{fmcf}$ we calculate
  \begin{align*}
    \lambda_{fmcf}
    &=\int_{\mathbb{T}^n}\nabla g(\nabla u)\cdot\nabla(Q\operatorname{div}_f v+\xi u)d\mu=\int_{\mathbb{T}^n}-\operatorname{div}_f(\nabla g(\nabla u))Q\operatorname{div}_fvd\mu+\xi\int_{\T^n}\n g(\n u)\cdot \n ud\mu\\
    &=\int_{\mathbb{T}^n}-D^2g(\nabla u):D^2u Q\operatorname{div}_fv+\int_{\mathbb{T}^n}\nabla f\cdot\nabla g(\nabla u)Q\operatorname{div}_fvd\mu+\xi\int_{\T^n}\n g(\n u)\cdot \n ud\mu\\
    &=\int_{\mathbb{T}^n}-g(\nabla u)|\operatorname{div}_fv|^2+\operatorname{div}_fv\left(g(\nabla u)\operatorname{div}_fv-QD^2g(\nabla u):D^2u\right)d\mu\\
    &\quad\quad+\int_{\mathbb{T}^n}\nabla f\cdot\nabla g(\nabla u)Q\operatorname{div}_fvd\mu+\xi\int_{\T^n}\n g(\n u)\cdot \n ud\mu\\
    &=\int_{\mathbb{T}^n}-g(\nabla u)|\operatorname{div}_fv|^2+\operatorname{div}_fv\operatorname{div}_f\left(g(\nabla u)v-Q\nabla g(\nabla u)\right)d\mu+\xi\int_{\T^n}\n g(\n u)\cdot \n ud\mu\\
    &=\int_{\mathbb{T}^n}-g(\nabla u)|\operatorname{div}_fv|^2+Dv^{\top}:D\left(g(\nabla u)v-Q\nabla g(\nabla u)\right)+vD^2f\left(g(\nabla u)v-Q\nabla g(\nabla u)\right)d\mu\\
    &\qquad +\xi\int_{\T^n}\n g(\n u)\cdot \n ud\mu\\
    &=\int_{\mathbb{T}^n}-g(\nabla u)|\operatorname{div}_fv|^2d\mu+vD^2f\left(g(\nabla u)v-Q\nabla g(\nabla u)\right)d\mu+\xi\int_{\T^n}\n g(\n u)\cdot \n ud\mu\\
    &\quad\quad+\int_{\mathbb{T}^n}Dv^{\top}:\left(v\otimes D^2u\nabla g(\nabla u)+g(\nabla u)Dv-\nabla g(\nabla u)\otimes D^2uv-QD^2g(\nabla u)D^2u \right)d\mu\\
    &=\int_{\mathbb{T}^n}-g(\nabla u)|\operatorname{div}_fv|^2d\mu+\int_{\mathbb{T}^n}vD^2f\left(g(\nabla u)v-Q\nabla g(\nabla u)\right)d\mu\\
    &\quad\quad+\int_{\mathbb{T}^n}Dv^{\top}:\left(g(\nabla u)Dv-QD^2g(\nabla u)D^2u\right)d\mu\\
    &\quad\quad+\int_{\mathbb{T}^n}Dv^{\top}:\left(v\otimes D^2u\nabla g(\nabla u)-\nabla g(\nabla u)\otimes D^2uv\right)d\mu+\xi\int_{\T^n}\n g(\n u)\cdot \n ud\mu\\
    &=\int_{\mathbb{T}^n}-g(\nabla u)|\operatorname{div}_fv|^2d\mu+\int_{\mathbb{T}^n}vD^2f\left(g(\nabla u)v-Q\nabla g(\nabla u)\right)d\mu\\
    &\quad\quad+\int_{\mathbb{T}^n}D^2u(\mathrm{Id}-v\otimes v):\left(\frac{g(\nabla u)}{Q^2}(\mathrm{Id}-v\otimes v)-D^2g(\nabla u)\right)D^2ud\mu\\
    &\quad\quad+\int_{\mathbb{T}^n}Dv^{\top}:\left(v\otimes D^2u\nabla g(\nabla u)-\nabla g(\nabla u)\otimes D^2uv\right)d\mu+\xi\int_{\T^n}\n g(\n u)\cdot \n ud\mu.
  \end{align*}
  For $\lambda_{pert}$ we calculate
  \begin{align*}
    \lambda_{pert}
    &=\int_{\mathbb{T}^n}\nabla g(\nabla u)\cdot\nabla\Delta_f u +\nabla Q\cdot D^2g(\nabla u)\nabla Q  d\mu\\
    &=\int_{\mathbb{T}^n}-D^2g(\nabla u)D^2u:D^2u+\nabla Q\cdot D^2g(\nabla u)\nabla Qd\mu-\int_{\mathbb{T}^n}\nabla g(\nabla u)D^2f\nabla ud\mu\\
    &=\int_{\mathbb{T}^n}-D^2u(\mathrm{Id}-v\otimes v):D^2g(\nabla u)D^2ud\mu-\int_{\mathbb{T}^n}\nabla g(\nabla u)D^2f\nabla ud\mu.
  \end{align*}
  For $\sigma $ we calculate
  \begin{align*}
    \sigma
    &=\int_{\mathbb{T}^n}\nabla g(\nabla u)\cdot\nabla Qd\mu=\int_{\mathbb{T}^n}\nabla g(\nabla u)\cdot D^2uvd\mu=\int_{\mathbb{T}^n}\nabla ( g(\nabla u))\cdot vd\mu\\
    &=-\int_{\mathbb{T}}g(\nabla u)\operatorname{div}_f vd\mu.
  \end{align*}
  Plug the above four terms into (\ref{Ito formula of I}) yields the result.
\end{proof}

\begin{remark}
  If we choose $g(p)$ as a function of $\sqrt{1+|p|^2}$, we will find that \[v\otimes D^2u\nabla g(\nabla u)-\nabla g(\nabla u)\otimes D^2uv=0.\]
  which will simplify the expression of $d\mathcal{I}$.
\end{remark}

Next we will choose $g$ as a function of $\Q$, which gives us more geometric meaning and is sufficiently to get the uniform estimates in $\epsilon$.

\begin{lemma}\label{lemma of expectation of I}
  Let $\epsilon\geq 0$, and $u$ be solution of (\ref{viscous equation}). Let $h\in C^2[1,\infty)$ be a non-negative, monotone
  and convex function with bounded second order derivative and $h^{\prime}(1)-h(1)\geq 0$ and 
  \[\mathcal{I}(t)=\int_{\mathbb{T}^n}h(\Q)d\mu,\quad\in I.\]
  For $q\in [1,2]$ we have 
  \begin{align*}
    \mathbb{E}&\mathcal{I}^q(t)\\ 
    &+\epsilon q\mathbb{E}\int^t_0\mathcal{I}^{q-1}(s)\int_{\mathbb{T}^n}\left(h^{\prime\prime}(Q(s))|D^2uv|^2(s)+\frac{h^{\prime}(Q(s))}{Q(s)}(|D^2u|^2-|D^2uv|^2(s))\right)d\mu ds\\
    &+\frac{2q-q^2}{2}\mathbb{E}\int^t_0\mathcal{I}^{q-1}(s)\int_{\mathbb{T}^n}h(Q(s))|\operatorname{div}_fv|^2(s)d\mu ds\\
    &+q\E\int^t_0\mathcal{I}^{q-1}\int_{\mathbb{T}^n}\nabla u\cdot\left[\left(\epsilon+\frac{h^{\prime}(Q(s))}{Q(s)}-\frac{h(Q(s))}{2Q(s)^2}\right)D^2f-\xi\frac{h^{\prime}(Q(s))}{Q(s)}\mathrm{Id}\right] \nabla ud\mu ds\\
    &+q\mathbb{E}\int^t_0\mathcal{I}^{q-1}(s)\int_{\mathbb{T}^n}\left(\frac{h^{\prime}(Q(s))}{Q(s)}-\frac{h(Q(s))}{2Q^2(s)}\right)\left(|D^2u|^2-2|D^2uv|^2+|v\cdot D^2u v|^2\right)(s)d\mu ds\\ 
    &+q\mathbb{E}\int^t_0\mathcal{I}^{q-1}(s)\int^t_0\int_{\mathbb{T}^n}h^{\prime\prime}(Q(s))\left(|D^2uv|^2-|v\cdot D^2uv|^2\right)(s)d\mu ds\\ 
    \leq&\mathbb{E}\mathcal{I}^q(0)\quad \forall t\in I.
  \end{align*}
\end{lemma}

\begin{lemma}\label{eigenvalue}
  Let $v$ be a $n$-dimensional vector. Then the eigenvalues of matrix $a \mathrm{Id}+b v\otimes v$ are $a$ and $a+b|v|^2$.
\end{lemma}
\begin{lemma}[\cite{Dabrock2021}*{Lemma 5.6}]\label{positive definite lemma}
  Let $A,B,C \in\mathbb{R}^{n\times n}$ be symmetric matrices with $B,C\geq 0$, Then 
  \[AB:CA=\sum_{ij}(AB)_{ij}(CA)_{ij}=\operatorname{tr}(AB)(CA)^{\top}\geq 0.\]
\end{lemma}

\begin{proof}[Proof of lemma \ref{lemma of expectation of I}]
  Let $g(p)=h(\mathbf{Q}(p))$ for $p\in\mathbb{R}^n$. Then 
  \[
    \begin{aligned}
      &\nabla g(p)=h^{\prime}(\mathbf{Q}(p))\mathbf{v}(p).\\
      &D^2 g(p)=h^{\prime \prime}(\mathbf{Q}(p))\mathbf{v}(p)\otimes\mathbf{v}(p)+\frac{h^{\prime}(\mathbf{Q}(p))}{\mathbf{Q}(p)}(\mathrm{Id}-\mathbf{v}(p)\otimes\mathbf{v}(p)).
    \end{aligned}
  \]
  At this time, two terms in the expression of $d\mathcal{I}$ in (\ref{Ito formula of I}) gets easier:
  \[
    \mathbf{v}(p)\otimes D^2u\nabla g(p)-\nabla g(p)\otimes D^2u\mathbf{v}(p)=0.
  \]
  \[
    \frac{1}{2}\mathbf{v}(p)D^2f\left(g(p)\mathbf{v}(p)-\mathbf{Q}(p)\nabla g(p)\right)=\frac{1}{2}(h(\mathbf{Q}(p))-\mathbf{Q}(p)h^{\prime}(\mathbf{Q}(p)))\mathbf{v}(p)\cdot D^2f\mathbf{v}(p).
  \]
  Moreover, 
  \begin{equation}\label{positive definite}
    \begin{aligned}
      &D^2g(p)-\frac{g(p)}{2Q^2}(\mathrm{Id}-\mathbf{v}(p)\otimes \mathbf{v}(p))\\
    =&\left(\frac{h^{\prime}(\mathbf{Q}(p))}{\mathbf{Q}(p)}-\frac{h(\mathbf{Q}(p))}{2\mathbf{Q}(p)^2}\right)\mathrm{Id}+\left(h^{\prime\prime}(\mathbf{Q}(p))-\frac{h^{\prime}(\mathbf{Q}(p))}{\mathbf{Q}(p)}+\frac{h(\mathbf{Q}(p))}{2\mathbf{Q}(p)^2}\right)\mathbf{v}(p)\otimes \mathbf{v}(p).
    \end{aligned}
  \end{equation}
  Note that \[\frac{d}{dx}(xh^{\prime}(x)-h(x))=xh^{\prime\prime}(x)\geq 0\quad \forall x\in[1,\infty).\] Thus function $xh^{\prime}(x)-h(x)$ is increasing with initial date $h^{\prime}(1)-h(1)\geq 0$, then 
  \[
    \frac{h^{\prime}(x)}{x}-\frac{h(x)}{x^2}\geq 0,\quad \forall x\geq 1.
  \]
  Now with lemma \ref{eigenvalue} we know that eigenvalues of (\ref{positive definite}) are given by:
  $$
    \begin{aligned}
      &\frac{h^{\prime}(Q(p))}{Q(p)}-\frac{h(Q(p))}{2Q(p^2}\geq 0\quad\text{and} \\ 
      &\frac{h^{\prime}{Q(p)}}{Q(p)^3}-\frac{h(Q(p))}{2Q(p)^4}+h^{\prime\prime}(Q(p))\frac{|p|^2}{Q(p)^2}\geq 0.
    \end{aligned}
  $$
  Which shows that (\ref{positive definite}) is non-negative. We adopt the same notation as lemma \ref{lemma of Ito formula of I}: $Q=\Q, v=\mathbf{v}(\nabla u)$.
  We can apply lemma \ref{lemma of Ito formula of I} to $\mathcal{I}(t)$ and deduce
  \[
    \begin{aligned}
      d\mathcal{I}
      &=-\int_{\mathbb{T}^n}\epsilon\left(h^{\prime\prime}(Q)|D^2uv|^2+\frac{h^{\prime}(Q)}{Q}(|D^2u|^2-|D^2uv|^2)\right)d\mu\\ 
      &\quad-\int_{\mathbb{T}^n}\frac{1}{2}h(Q)|\operatorname{div}_fv|^2d\mu-\int_{\mathbb{T}^n}\left(\epsilon+\frac{h^{\prime}(Q)}{Q}-\frac{h(Q)}{2Q^2}\right)\nabla u\cdot D^2f\nabla ud\mu\\ 
      &\quad -\int_{\mathbb{T}^n}D^2u(\mathrm{Id}-v\otimes v):\left(\left(\frac{h^{\prime}(Q)}{Q}-\frac{h(Q)}{2Q^2}\right)(\mathrm{Id}-v\otimes v)+h^{\prime\prime}(Q) v\otimes v\right)D^2ud\mu dt\\ 
      &\quad+\xi\int_{\T^n}\frac{h^{\prime}(Q)}{Q}|\n u|^2d\mu-\int_{\mathbb{T}^n}h(Q)\operatorname{div}_f vd\mu dW\\
      &=-\int_{\mathbb{T}^n}\epsilon\left(h^{\prime\prime}(Q)|D^2uv|^2+\frac{h^{\prime}(Q)}{Q}(|D^2u|^2-|D^2uv|^2)\right)d\mu\\ 
      &\quad-\int_{\mathbb{T}^n}\frac{1}{2}h(Q)|\operatorname{div}_fv|^2d\mu-\int_{\mathbb{T}^n}\nabla u\cdot\left[\left(\epsilon+\frac{h^{\prime}(Q)}{Q}-\frac{h(Q)}{2Q^2}\right)D^2f-\xi\frac{h^{\prime}(Q)}{Q}\mathrm{Id}\right] \nabla ud\mu\\ 
      &\quad -\int_{\mathbb{T}^n}D^2u(\mathrm{Id}-v\otimes v):\left(\left(\frac{h^{\prime}(Q)}{Q}-\frac{h(Q)}{2Q^2}\right)(\mathrm{Id}-v\otimes v)+h^{\prime\prime}(Q) v\otimes v\right)D^2ud\mu dt\\ 
      &-\int_{\mathbb{T}^n}h(Q)\operatorname{div}_f vd\mu dW.
    \end{aligned}
  \]
By lemma \ref{positive definite lemma}  the forth integral in the previous formula is non-negative. We can apply Fatou's lemma to get rid of the locality and deduce that 
\[
    \begin{aligned}
      \mathbb{E}\int_{\mathbb{T}^n}&h(Q(t))d\mu\\
      &+\epsilon\mathbb{E}\int^t_0\int_{\mathbb{T}^n}\left(h^{\prime\prime}(Q(s))|D^2uv|^2(s)+\frac{h^{\prime}(Q(s))}{Q(s)}\left(|D^2u|^2(s)-|D^2uv|^2(s)\right)\right)d\mu ds\\ 
      &+\frac{1}{2}\mathbb{E}\int^t_0\int_{\mathbb{T}^n}h(Q(s))|\operatorname{div}_fv|^2(s)d\mu ds\\
      &+\E\int^t_0\int_{\mathbb{T}^n}\nabla u\cdot\left[\left(\epsilon+\frac{h^{\prime}(Q)}{Q}-\frac{h(Q)}{2Q^2}\right)D^2f-\xi\frac{h^{\prime}(Q)}{Q}\mathrm{Id}\right] \nabla ud\mu ds\\ 
      &+\mathbb{E}\int^t_0\int_{\mathbb{T}^n}\left(\frac{h^{\prime}(Q(s))}{Q(s)}-\frac{h(Q(s))}{2Q^2(s)}\right)\left(|D^2u|^2-2|D^2uv|^2+|v\cdot D^2u v|^2\right)(s)d\mu ds\\ 
      &+\mathbb{E}\int^t_0\int_{\mathbb{T}^n}h^{\prime\prime}(Q(s))\left(|D^2uv|^2-|v\cdot D^2uv|^2\right)(s)d\mu ds\\ 
      \leq \mathbb{E}&\int^t_0\int_{\mathbb{T}^n}h(Q(0))d\mu\quad \forall t\in I.
    \end{aligned}
\]
Where \[\left(|D^2u|^2-2|D^2uv|^2+|v\cdot D^2u v|^2\right)=D^2u(\mathrm{Id}-v\otimes v):(\mathrm{Id}-v\otimes v)D^2u\geq 0.\]
For $q\in[1,2]$ we want to use the It\^{o} formula for the function $x\mapsto |x|^q.$ This function is not twice continuously differentiable for $q<2$, 
so the classical It\^{o} formula does not apply directly. Nevertheless, we can first do the calculation  to the function
$x\mapsto (x+\theta)^q$ for $\theta>0$ which is twice continuously differentiable on $[0,\infty)$, then send $\theta\rightarrow 0$. 
We infer
$$
    \begin{aligned}
      d\mathcal{I}^q=
      &-\epsilon q\mathcal{I}^{q-1}\int_{\mathbb{T}^n}\left(h^{\prime\prime}(Q)|D^2uv|^2+\frac{h^{\prime}(Q)}{Q}(|D^2u|^2-|D^2uv|^2)\right)d\mu\\
      &-\frac{q}{2}\mathcal{I}^{q-1}\int_{\mathbb{T}^n}h(Q)|\operatorname{div}_fv|^2d\mu \\
      &-q\mathcal{I}^{q-1}\int_{\mathbb{T}^n}\nabla u\cdot\left[\left(\epsilon+\frac{h^{\prime}(Q)}{Q}-\frac{h(Q)}{2Q^2}\right)D^2f-\xi\frac{h^{\prime}(Q)}{Q}\mathrm{Id}\right] \nabla ud\mu\\
      &-q\mathcal{I}^{q-1}\int_{\mathbb{T}^n}D^2u(\mathrm{Id}-v\otimes v):\left(\left(\frac{h^{\prime}(Q)}{Q}-\frac{h(Q)}{2Q^2}\right)(\mathrm{Id}-v\otimes v)+h^{\prime\prime}(Q) v\otimes v\right)D^2ud\mu \\
      &+\frac{q(q-1)}{2}\mathcal{I}^{q-2}\left(\int_{\mathbb{T}^n}h(Q)\operatorname{div}_fvd\mu\right)^2dt\\
      &-q\mathcal{I}^{q-1} \int_{\mathbb{T}^n}h(Q)\operatorname{div}_f vd\mu dW.\\
    \leq &-\epsilon q\mathcal{I}^{q-1}\int_{\mathbb{T}^n}\left(h^{\prime\prime}(Q)|D^2uv|^2+\frac{h^{\prime}(Q)}{Q}(|D^2u|^2-|D^2uv|^2)\right)d\mu\\
      &+\left(-\frac{q}{2}+\frac{q(q-1)}{2}\right)\int_{\mathbb{T}^n}h(Q)|\operatorname{div}_fv|^2d\mu\\ 
      &-q\mathcal{I}^{q-1}\int_{\mathbb{T}^n}\nabla u\cdot\left[\left(\epsilon+\frac{h^{\prime}(Q)}{Q}-\frac{h(Q)}{2Q^2}\right)D^2f-\xi\frac{h^{\prime}(Q)}{Q}\mathrm{Id}\right] \nabla ud\mu\\
      &-q\mathcal{I}^{q-1}\int_{\mathbb{T}^n}D^2u(\mathrm{Id}-v\otimes v):\left(\left(\frac{h^{\prime}(Q)}{Q}-\frac{h(Q)}{2Q^2}\right)(\mathrm{Id}-v\otimes v)+h^{\prime\prime}(Q) v\otimes v\right)D^2ud\mu dt \\
      &-q\mathcal{I}^{q-1} \int_{\mathbb{T}^n}h(Q)\operatorname{div}_f vd\mu dW. 
    \end{aligned}
$$
As before, since the integral defines a local martingale and using Fatou's lemma, we get 
\begin{align*}
    \mathbb{E}&\mathcal{I}^q(t)\\ 
    &+\epsilon q\mathbb{E}\int^t_0\mathcal{I}^{q-1}(s)\int_{\mathbb{T}^n}\left(h^{\prime\prime}(Q(s))|D^2uv|^2(s)+\frac{h^{\prime}(Q(s))}{Q(s)}(|D^2u|^2-|D^2uv|^2(s))\right)d\mu ds\\
    &+\frac{2q-q^2}{2}\mathbb{E}\int^t_0\mathcal{I}^{q-1}(s)\int_{\mathbb{T}^n}h(Q(s))|\operatorname{div}_fv|^2(s)d\mu ds\\
    &+q\E\int^t_0\mathcal{I}^{q-1}\int_{\mathbb{T}^n}\nabla u\cdot\left[\left(\epsilon+\frac{h^{\prime}(Q(s))}{Q(s)}-\frac{h(Q(s))}{2Q(s)^2}\right)D^2f-\xi\frac{h^{\prime}(Q(s))}{Q(s)}\mathrm{Id}\right] \nabla ud\mu ds\\
    &+q\mathbb{E}\int^t_0\mathcal{I}^{q-1}(s)\int_{\mathbb{T}^n}\left(\frac{h^{\prime}(Q(s))}{Q(s)}-\frac{h(Q(s))}{2Q^2(s)}\right)\left(|D^2u|^2-2|D^2uv|^2+|v\cdot D^2u v|^2\right)(s)d\mu ds\\ 
    &+q\mathbb{E}\int^t_0\mathcal{I}^{q-1}(s)\int^t_0\int_{\mathbb{T}^n}h^{\prime\prime}(Q(s))\left(|D^2uv|^2-|v\cdot D^2uv|^2\right)(s)d\mu ds\\ 
    \leq&\mathbb{E}\mathcal{I}^q(0)\quad \forall t\in I.
\end{align*}
\end{proof}

\begin{proof}[Proof of proposition \ref{uniform estimate of energy}]
  We choose $h(x)=\alpha\epsilon x^2$ with $\alpha$ small enough such that  $0<\dfrac{4\alpha}{2+3\alpha}\leq \delta$ holds in Lemma \ref{lemma of expectation of I}. Then $h(\mathbf{Q}(p))=\alpha\epsilon(1+|p|^2)$, note that $d\int_{\T^n}|\n u|^2d\mu=d\int_{\T^n}g(\Q)\mu$. Then for $q=1$,
  we get
  \[
    \begin{aligned}
      \mathbb{E}\|\nabla u\|^2_{L^2(\mathbb{T}^n,d\mu)}&+2\epsilon\mathbb{E}\int^t_0\int_{\mathbb{T}^n}|D^2u(s)|^2d\mu ds+\frac{1}{2}\mathbb{E}\int^t_0\int_{\mathbb{T}^n}\mathbf{Q}(\n u(s))^2|\operatorname{div}_f v(s)|^2d\mu ds\\ 
      &\E\int^t_0\int_{\T^n}\n u\cdot\left(\left(1+\frac{3\alpha}{2}\right)D^2f-2\alpha\xi\mathrm{Id}\right)\n ud\mu ds\\ 
      &+\mathbb{E}\int^t_0\int_{\mathbb{T}^n}\left(\frac{3}{2}|D^2u|^2-|D^2uv|^2-\frac{1}{2}|v\cdot D^2u v|^2\right)(s)d\mu ds\\ 
      \leq &\mathbb{E}\|\nabla u_0\|_{L^2(\mathbb{T}^n,d\mu)}.
    \end{aligned}
  \]
   Notice that all the terms on the left-hand side in Lemma \ref{lemma of expectation of I} are non-negative under the condition $D^2f\geq \delta\xi\mathrm{Id}$ if we choose $h(x)=\alpha\epsilon x^2$ with $\alpha$ satisfying $\dfrac{4\alpha}{2+3\alpha}\leq \delta$. Therefore, we obtain in particular a uniform bound with respect to $t$ of \[\mathbb{E}\int^t_0\mathcal{I}^{q-1}(s)\int_{\mathbb{T}^n}h(Q(s))|\operatorname{div}_fv|^2(s)d\mu ds,\]  which we will use in the following calculation.
  For the stochastic integral we apply the Burkholder-Davis-Gundy inequality and obtain
\[
    \begin{aligned}
      \mathbb{E}&\sup_{t\in I}\left[\int^t_0\mathcal{I}^{q-1}(s)\int_{\mathbb{T}^n}h(Q(s))\operatorname{div}_fv(s)d\mu dW(s)\right]\\ 
      &\leq C\mathbb{E}\left[\int^{\sup I}_0\left(\mathcal{I}^{q-1}(s)\int_{\mathbb{T}^n}h(Q(s))\operatorname{div}_fv(s)d\mu\right)^2dt\right]^{1/2}\\ 
      &\leq C\mathbb{E}\left[\int^{\sup I}_0\left(\mathcal{I}^{2q-1}(s)\int_{\mathbb{T}^n}h(Q(s))|\operatorname{div}_fv(s)|^2d\mu\right)dt\right]^{1/2}\\ 
      &\leq C\mathbb{E}\left[\sup_{t\in I}\mathcal{I}^{q}(t)\int^{\sup I}_0\left(\mathcal{I}^{q-1}(s)\int_{\mathbb{T}^n}h(Q(s))|\operatorname{div}_fv(s)|^2d\mu\right)dt\right]^{1/2}\\
      &\leq \frac{C\delta}{2}\mathbb{E}\sup_{t\in I}\mathcal{I}^q(t)+\frac{C}{2\delta}\mathbb{E}\sup_{t\in I}\mathcal{I}^{q}(t)\int^{\sup I}_0\left(\mathcal{I}^{q-1}(s)\int_{\mathbb{T}^n}h(Q(s))|\operatorname{div}_fv(s)|^2d\mu\right)dt\\ 
      &\leq \frac{C\delta}{2}\mathbb{E}\sup_{t\in I}\mathcal{I}^q(t)+\frac{C}{\delta(2q-q^2)}\mathbb{E}\mathcal{I}^q(0).
    \end{aligned}
\]
for any $\delta >0.$ Thus, 
\[\mathbb{E}\sup_{t\in I}\mathcal{I}^q(t)\leq \mathbb{E}\mathcal{I}^q(0)+\frac{C\delta}{2}\mathbb{E}\sup_{t\in I}\mathcal{I}^q(t)+\frac{C}{\delta(2q-q^2)}\mathbb{E}\mathcal{I}^q(0).\]
Now, we choose $\delta=\frac{1}{C}$, Then we finally get
\[\mathbb{E}\sup_{t\in I}\mathcal{I}^q(t)\leq \left(2+\frac{2C^2}{2q-q^2}\right)\mathbb{E}\mathcal{I}^q(0).\]
That is 
  \[\mathbb{E}\sup_{t\in I}\|\nabla u\|^{2q}_{L^2(\T^n,d\mu)}\leq\left(2+\frac{2C^2}{2q-q^2}\right)\mathbb{E}\|\nabla u_0\|^{2q}_{L^2(\T^n,d\mu)}.\]
  Now let $\mathbb{P}$ - esssup $\|\nabla u\|_{L^{\infty}\left(I ; L^{\infty}\left(\mathbb{T}^n\right)\right)}=L<\infty$. Then we can estimate
\[
|\mathbf{v}(\nabla u)|=\frac{|\nabla u|}{\mathbf{Q}(\nabla u)} \leq \frac{L}{\sqrt{1+L^2}}<1.
\]
and

\begin{align*}
3\left|D^2 u\right|^2 & -2\left|D^2 u \mathbf{v}(\nabla u)\right|^2-\left|\mathbf{v}(\nabla u) \cdot D^2 u \mathbf{v}(\nabla u)\right|^2 \\
& \geq\left(3-2|\mathbf{v}(\nabla u)|^2-|\mathbf{v}(\nabla u)|^4\right)\left|D^2 u\right|^2 \\
& \geq \frac{3\left(1+L^2\right)^2-2 L^2\left(1+L^2\right)-L^4}{\left(1+L^2\right)^2}\left|D^2 u\right|^2 \\
& =\frac{3+4 L^2}{\left(1+L^2\right)^2}\left|D^2 u\right|^2 .
\end{align*}
Hence,
\begin{align*}
  \mathbb{E}\|\nabla u(t)\|_{L^2\left(\mathbb{T}^n,d\mu\right)}^2 & +\frac{3+4 L^2}{2\left(1+L^2\right)^2} \mathbb{E} \int_0^t \int_{\mathbb{T}^n}\left|D^2 u(s)\right|^2 \mathrm{~d} s \\
  & \leq \mathbb{E}\left\|\nabla u_0\right\|_{L^2\left(\mathbb{T}^n,d\mu\right)}^2 \forall t \in I .
\end{align*}

\end{proof}

\begin{proof}[Proof of proposition \ref{Maximum principle}]
  For $M>0$, Let $h_M\in C^2([1,\infty])$ be a modification of $x\mapsto (x-M)_+$ with the following properties:
  \begin{itemize}
    \item $h_M\geq 0,$
    \item $h_M$ is monotone increasing and convex,
    \item $h^{\prime\prime}_M$ is bounded,
    \item $h^{\prime}_M(1)-h_M(1)\geq 0$ and
    \item $h_M(x)>0\Leftrightarrow x>M.$
  \end{itemize}
  For example, one could choose
  \[
    h_M(x):= 
    \begin{cases}
      0 & x \leq M, \\ 
      \alpha\epsilon(x-M)^2 & x \in(M, M+1), \\ 
      \alpha\epsilon(2 x-2 M-1) & x \in[M+1, \infty).
    \end{cases}
  \]
  with $\alpha$ satisfying $0<\dfrac{4\alpha}{2+3\alpha}\leq \delta$.
  From lemma \ref{lemma of expectation of I} we deduce
  $$\mathbb{E}\int_{\mathbb{T}^n}h_M(Q(\nabla u(t)))d\mu\leq \mathbb{E}\int_{\mathbb{T}^n}h_M(Q(\nabla u_0))d\mu\quad \forall t\in I.$$
  suppose $\|\nabla u_0\|_{L^{\infty}(\mathbb{T}^n)}\leq L$ $\mathbb{P}$-a.s. then we have 
  $$\|Q(\nabla u_0)\|_{L^{\infty}(\mathbb{T}^n)}\leq\sqrt[]{1+L^2}\quad \mathbb{P}-a.s.$$
  Hence 
  $$\mathbb{E}\int_{\mathbb{T}^n}h_{\sqrt[]{1+L^2}}(Q(\nabla u(t)))d\mu=0\quad \forall t\in I$$
  which implies $\|Q(\nabla u(t))\|_{L^{\infty}(\mathbb{T}^n)}\leq \sqrt[]{1+L^2}\quad \mathbb{P}$-a.s. for all $t\in I.$
  Therefore,
  $$\|\nabla u\|_{L^{\infty}(I;L^{\infty}(\mathbb{T}^n))}\leq L.$$

\end{proof}
\begin{remark}
  We can also choose $h(x)=\epsilon x$ in Lemma \ref{lemma of expectation of I} to deduce the $q$-th moment of the area. If we choose $q=0$, we get 
  \begin{align*}
    \E\int_{\T^n}&\mathbf{Q}(\n u(t))d\mu+\epsilon\E\int^t_0\int_{\T^n}\frac{1}{\Q}\left(|D^2u|^2-|D^2u\V|^2\right)(s)d\mu ds\\
    &+\frac{1}{2}\E\int^t_0\int_{\T^n}\Q|\mathrm{div}_f\V|^2(s)d\mu ds\\
    &+\E\int^t_0\int_{\T^n}\n u\cdot\left(1+\frac{1}{\Q}D^2f-\frac{1}{\Q}\xi\mathrm{Id}\right)\n u(s)d\mu ds\\
    &+\E\int^t_0\int_{\T^n}\frac{1}{2\Q}\left(|D^2u|^2-2|D^2u\V|^2+|\V\cdot D^2u\V|^2\right)(s)d\mu ds\\
    &\leq \E\int_{\T^n}\mathbf{Q}(\n u_0)d\mu.
  \end{align*}
\end{remark}

\section{Existence of martingale solutions of Stochastic $f$-MCF}\label{sec:6}
We now set about establishing the existence of a martingale solution of the stochastic $f$-mean curvature flow (\ref{SfMCF}) by passing the limit $\epsilon\rightarrow 0$ with the help of the above uniform estimates. The proof is standard and similar with those in \cites{Dabrock2020,Hofmanova2017}.
\begin{proof}[Proof of Theorem \ref{main theorem}]
  By Theorem \ref{existence of viscous martingale solution} we deduce that for $\epsilon >0$ there exists a martinagle solution $u^{\epsilon}$ of (\ref{viscous equation}) with initial data $\Lambda$. Since it is Theorem \ref{existence} we use to get the existence of $u^{\epsilon}$ which is proved by Jakowbuski-Skorohkod representation theorem \cite{Jakubowski1997}, we can fix  one probability space $(\Omega,\mathcal{F},\p)=([0,1],\mathcal{B}([0,1]),\mathcal{L})$ such that for each $\epsilon>0$ we can find
  \begin{itemize}
    \item a normal filtration $(\mathcal{F}^{\epsilon}_t)_{t\in[0,\infty)}$,
    \item a real-valued $(\mathcal{F}^{\epsilon}_t)$-Wiener process $W^{\epsilon}$ and
    \item a $(\mathcal{F}^{\epsilon}_t)$-predictable process $u^{\epsilon}$ with $u^{\epsilon}\in L^2(\Omega;L^2(0,T;H^2(\T^n,d\mu)))$ for all $T\in [0,\infty)$.
  \end{itemize}
  such that 
  \begin{equation}\label{integral}
    \begin{aligned}
      u^{\epsilon}(t)-u^{\epsilon}(0)&=\int^t_0\epsilon\Delta_f u^{\epsilon}(s)+\mathbf{Q}(\n u^{\epsilon}(s))\mathrm{div}_f\mathbf{v}(\n u^{\epsilon}(s))+\xi u^{\epsilon}(s)ds\\
      &\quad+\int^t_0\mathbf{Q}(\n u^{\epsilon}(s))dW^{\epsilon}_s\quad \text{in}\,L^2(\T^n,d\mu)\, \forall t\in\: [0,\infty).
    \end{aligned}
  \end{equation}
  and $\p\circ (\n u^{\epsilon}(0))^{-1}=\Lambda$. With the assumption on the support of $\Lambda$ we have $\|\n u^{\epsilon}(0)\|_{L^{\infty}}\leq L\,\p$-a.s. By Proposition \ref{Maximum principle} we have 
  \[
    \|\n u^{\epsilon}\|L^{\infty}(0,\infty;L^{\infty})\leq L\,\;\p\text{-a.s.}.
  \]
  By Proposition \ref{uniform estimate of energy}, we have for $q\in [1,2)$
  \begin{equation}\label{uniform upper bound of D^2u}
    \begin{aligned}
      & \left\|\nabla u^{\epsilon}\right\|_{L^{2 q}\left(\Omega ; C\left([0, \infty) ; L^2\left(\mathbb{T}^n,d\mu\right)\right)\right)}^{2 q} \leq\left(2+\frac{2 C^2}{2 q-q^2}\right) \mathbb{E}\left\|\nabla u^{\epsilon}(0)\right\|_{L^2\left(\mathbb{T}^n,d\mu\right)}^{2 q} \leq C_{q, L} , \\
      & \left\|\mathrm{D}^2 u^{\epsilon}\right\|_{L^2\left(\Omega ; L^2\left(0, \infty ; L^2\left(\mathbb{T}^n,d\mu\right)\right)\right)} \leq C_L\left\|u^{\epsilon}(0)\right\|_{L^2\left(\Omega ; H^1\left(\mathbb{T}^n,d\mu\right)\right)} .
    \end{aligned}
  \end{equation}
  By It\^{o} formula \ref{Ito formula}, we have 
  \begin{align*}
    d\|u^{\epsilon}\|^2_{L^2(\T^n)}&=\int_{\T^n}2u^{\epsilon}\cdot\left(\epsilon\Delta_f u^{\epsilon}+\mathbf{Q}(\n u^{\epsilon})\mathrm{div}_f\mathbf{v}(\n u^{\epsilon})+\frac{1}{2}\mathbf{v}(\n u^{\epsilon})\cdot D^2u\mathbf{v}(\n u^{\epsilon})+\xi u^{\epsilon}\right)d\mu dt\\
    &\quad+\int_{\T^n}\mathbf{Q}(\n u^{\epsilon})^2d\mu dt +2\int_{\T^n}u^{\epsilon}\mathbf{Q}(\n u^{\epsilon})d\mu dW\\
    &=\int_{\T^n}2u^{\epsilon}\left((1+\epsilon)\Delta_fu^{\epsilon}-\frac{1}{2}\mathbf{v}(\n u^{\epsilon})\cdot D^2u\mathbf{v}(\n u^{\epsilon})+\xi u^{\epsilon}\right)d\mu dt\\
    &\quad+\int_{\T^n}\mathbf{Q}(\n u^{\epsilon})^2d\mu dt +2\int_{\T^n}u^{\epsilon}\mathbf{Q}(\n u^{\epsilon})d\mu dW\\
    &\leq \left(\frac{1}{2}+2\xi\right)\|u^{\epsilon}\|^2_{L^2(\T^n,d\mu)}dt+\frac{1}{2}\|D^2u\|^2_{L^2(\T^n,d\mu)}+2\int_{\T^n}u^{\epsilon}\mathbf{Q}(\n u^{\epsilon})d\mu dW.
  \end{align*}
  By the Burkholder-Davis-Gundy inequality, we can get the upper bound of the supremum of $\|u^{\epsilon}\|^2_{L^2(\T^n,d\mu)}$. That is,
  \begin{align*}
    \E&\sup_{s\in[0,t]}\|u^{\epsilon}(s)\|^2_{L^2(\T^n,d\mu)}-\E\|u^{\epsilon}\|^2_{L^2(\T^n,d\mu)}\\
    &\leq \left(\frac{1}{2}+2\xi\right)\E\int^t_0\|u^{\epsilon}(s)\|^2_{L^2(\T^n,d\mu)}ds+\frac{1}{2}\E\int^t_0\|D^2u\|^2_{L^2(\T^n,d\mu)}ds\\
    &\qquad+C\E\left[\int^t_0\|u^{\epsilon}(s)\|^2_{L^2(\T^n,d\mu)}\|\mathbf{Q}(\n u^{\epsilon}(s))\|^2_{L^2(\T^n,d\mu)}ds\right]^{1/2}\\
    &\leq \left(\frac{1}{2}+2\xi\right)\E\int^t_0\|u^{\epsilon}(s)\|^2_{L^2(\T^n,d\mu)}ds+\frac{1}{2}\E\int^t_0\|D^2u\|^2_{L^2(\T^n,d\mu)}ds\\
    &\qquad +\frac{1}{2}\E\sup_{s\in[0,t]}\|u^{\epsilon}(s)\|^2_{L^2(\T^n,d\mu)}+C\E\int^t_0\|\mathbf{Q}(\n u^{\epsilon}(s))\|^2_{L^2(\T^n,d\mu)}ds.
  \end{align*}
  By (\ref{uniform upper bound of D^2u}) we have for all $T>0$ and $t\in [0,T]$,
  \[
    \E\sup_{s\in[0,t]}\|u^{\epsilon}(s)\|^2_{L^2(\T^n,d\mu)}\leq C\E\int^t_0\|u^{\epsilon}(s)\|^2_{L^2(\T^n,d\mu)}ds+C,
  \]
  with constant $C$ depends only on initial data $\Lambda$ and $T,\xi$. By Gronwall inequatlity, we have 
  \[
    \E\sup_{t\in[0,T]}\|u^{\epsilon}(t)\|^2_{L^2(\T^n,d\mu)}\leq C(\Lambda,T,\xi).
  \]
  With (\ref{uniform upper bound of D^2u}) in hand, the usual integral in (\ref{integral}) in bounded in $L^2(\Omega,;C^{0,1/2}([0,T];L^2(\T^n,d\mu)))$. By the factorization method \cite{Seidler1993}*{Theorem 1.1} and (\ref{integral}), we know that there exists a $\lambda>0$ such that the stochastic integral in (\ref{integral}) is unformly bounder in $L^2(\Omega;C^{0,\lambda}([0,T];L^2(\T^n,d\mu)))$. Thus, for some $\lambda\in (0,1/2)$, we have 
  \[
    \E\|u^{\epsilon}\|_{C^{0,\lambda}([0,T];L^2(\T^n,d\mu))}\leq C(\Lambda,T,\xi),
  \]
  unifromly in $\epsilon$. Taken the above together, we have $(u^{\epsilon})_{\epsilon>0}$ is uniformly bounded in 
  \[
    L^2\left(\Omega ; L^2\left([0, T] ; H^2\left(\mathbb{T}^n,d\mu\right)\right) \cap C\left([0, T] ; H^1\left(\mathbb{T}^n,d\mu\right)\right) \cap C^{0, \lambda}\left([0, T] ; L^2\left(\mathbb{T}^n,d\mu\right)\right)\right).
  \]
   and $(W^{\epsilon})_{\epsilon>0}$ is uniformly bounded in $L^2(\Omega;C^{0,\lambda}([0,T];\R))$.

   Next, we show that those uniform bounds imply the convergence of a subsequence in weak sence, then we identify the limit with the solution of (\ref{SfMCF}).

   We claim that the following two embeddings are compact
   \begin{align*}
    C^{0, \lambda}\left([0, T] ; L^2\left(\mathbb{T}^n,d\mu\right)\right)\cap C\left([0, T] ; H^1\left(\mathbb{T}^n,d\mu\right)\right)&\hookrightarrow C\left([0, T] ; (H^1\left(\mathbb{T}^n,d\mu\right),w)\right),\\
    L^2\left([0, T] ; H^2\left(\mathbb{T}^n,d\mu\right)\right)\cap C^{0, \lambda}\left([0, T] ; L^2\left(\mathbb{T}^n,d\mu\right)\right)&\hookrightarrow L^2\left([0, T] ; H^1\left(\mathbb{T}^n,d\mu\right)\right).
   \end{align*}
   The first embedding is compact due to Arzela-Ascoli theorem and the fact that the topology of $L^2\left(\mathbb{T}^n,d\mu\right)$ coincides with the weak topology $(H^1\left(\mathbb{T}^n,d\mu\right),w)$ of $H^1\left(\mathbb{T}^n,d\mu\right)$ on the bounded sets of $C\left([0, T] ; H^1\left(\mathbb{T}^n,d\mu\right)\right)$. The second embedding is compact due to the characterization of compacts sets in \cite{Simon}. 

   Thus, the joint laws of $(u^{\epsilon},W^{\epsilon})$ are tight in 
   \[
    \mathcal{X}^T_u\times\mathcal{X}^T_W,
   \]
   with
   \begin{align*}
    & \mathcal{X}_u^T:=C\left([0, T] ;\left(H^1\left(\mathbb{T}^n,d\mu\right), w\right)\right) \cap L^2\left(0, T ; H^1\left(\mathbb{T}^n,d\mu\right)\right) \cap\left(L^2\left(0, T ; H^2\left(\mathbb{T}^n,d\mu\right)\right), w\right), \\
    & \mathcal{X}_W^T:=C([0, T] ; \mathbb{R}).
  \end{align*}
  Since $T>0$ is arbitrary this also implies that the joint laws of $(u^{\epsilon},W^{\epsilon})$ are tight in $\mathcal{X}_u\times\mathcal{X}_W$ with
  \begin{align*}
    \mathcal{X}_u:= & C_{\mathrm{loc}}\left([0, \infty) ;\left(H^1\left(\mathbb{T}^n,d\mu\right), w\right)\right) \cap L_{\mathrm{loc}}^2\left(0, \infty ; H^1\left(\mathbb{T}^n,d\mu\right)\right)  \cap\left(L_{\mathrm{loc}}^2\left(0, \infty ; H^2\left(\mathbb{T}^n,d\mu\right)\right), w\right), \\
    \mathcal{X}_W:= & C_{\mathrm{loc}}([0, \infty) ; \mathbb{R}) .
  \end{align*}
  The Jakubowski-Skorohkod representation theorem for tight sequence in nonmetric spaces \cite{Jakubowski1997}*{Theorem 2} tells us that there exists a sequence $\epsilon_k\to 0$, a probability space $(\tilde{\Omega},\tilde{\mathcal{F}},\tilde{\mathbb{P}})$ and $\mathcal{X}_u\times\mathcal{X}_W$-valued random variables $(\tilde{u}^{k},\tilde{W}^{k})$ for $k\in\N$ and $(\tilde{u},\tilde{W})$ such that $\tilde{u}^k\to \tilde{u}$ a.s. in $\mathcal{X}_u$ and $\tilde{W}^k\to \tilde{W}$ a.s. in $\mathcal{X}_W$ and the joint law of $(\tilde{u}^{k},\tilde{W}^{k})$ coincide with the joint law of $(u^{\epsilon_k},W^{\epsilon_k})$ for $k\in\N$.
  Let 
  \begin{align*}
    & \tilde{\mathcal{F}}_t^k:=\bigcap_{s>t} \sigma\left(\left.\tilde{u}^k\right|_{[0, s]},\left.\tilde{W}^k\right|_{[0, s]},\{A \in \tilde{\mathcal{F}} \mid \tilde{\mathbb{P}}(A)=0\}\right), t \in[0, \infty), k \in \mathbb{N} \\
    & \tilde{\mathcal{F}}_t:=\bigcap_{s>t} \sigma\left(\left.\tilde{u}\right|_{[0, s]},\left.\tilde{W}\right|_{[0, s]},\{A \in \tilde{\mathcal{F}} \mid \tilde{\mathbb{P}}(A)=0\}\right), t \in[0, \infty) .
  \end{align*}
  We can prove that $\tilde{W}^k$ is a real-valued $(\tilde{\mathcal{F}}^k_t)_t$-Wiener process and $\tilde{u}^k$ is a solution of the viscous equation (\ref{viscous equation}) for $\epsilon_k$ and the Wiener process $\tilde{W}_k$.

  Let $t\in [0,\infty)$, for any $e\in H^2(\T^n,d\mu)$ we define 
  \begin{align*}
    \tilde{M}(t)&:=\left\langle \tilde{u}(t)-u_0,e\right\rangle_{H^1}-\int^t_0\left\langle\Delta_f \tilde{u}-\frac{1}{2}\mathbf{v}(\n \tilde{u})D^2\tilde{u}\mathbf{v}(\n \tilde{u})+\xi \tilde{u}(s),\:e\right\rangle_{H^{-2},H^2} ds,\\
    \tilde{M}^k(t)&:=\left\langle \tilde{u}^k(t)-u_0,e\right\rangle_{H^1}-\int^t_0\left\langle1+\epsilon_k)\Delta_f \tilde{u}^k-\frac{1}{2}\mathbf{v}(\n \tilde{u}^k)D^2\tilde{u}^k\mathbf{v}(\n \tilde{u}^k)+\xi \tilde{u}^k,\:e\right\rangle_{H^{-2},H^2}ds,\\
    M^k(t)&:=\left\langle u^k(t)-u_0,e\right\rangle_{H^1}-\int^t_0\left\langle(1+\epsilon_k)\Delta_f u^k-\frac{1}{2}\mathbf{v}(\n u^k)D^2u^k\mathbf{v}(\n u^k)+\xi u^k(s),e\right\rangle_{H^{-2},H^2}ds.
  \end{align*}
  For $t\in [0,\infty)$, we have 
  \[
    M^k(t)=\int^t_0\left\langle\mathbf{Q}(\n u^k) dW_s,\:e\right\rangle_{H^1}.
  \]
  For $s\in [0,t]$, we define 
  \[\gamma:C([0,s];(H^2(\T^n,d\mu),w))\times C([0,s],(\R,w))\to\R,\]
  be a bounded and continuous function. We use the following abbreviations 
  \[\gamma^k:=\gamma\left(u^k|_{[0,s]},W^k|_{[0,s]}\right),\quad\tilde{\gamma}^k:=\left(\tilde{u}^k|_{[0,s]},\tilde{W}^k|_{[0,s]}\right),\quad\tilde{\gamma}:=\gamma\left(\tilde{u}|_{[0,s]},\tilde{W}|_{[0,s]}\right).\]
  Since the joint law of $(\tilde{u}^k,\tilde{W}^k)$ coincides with the joint law of $(u^k,W^k)$, we have 
  \begin{equation}\label{28}
    \tilde{\E}\left(\left.\tilde{\gamma}^k\tilde{W}^k\right|^t_s\right)=0,\quad \tilde{\E}\left(\left.\tilde{\gamma}^k (\tilde{W}^k)^2\right|^t_s\right)=t-s.
  \end{equation}
Thus, we have
\begin{equation}\label{29}
  \begin{aligned}
    0&=\E\left(\left.\gamma^k M^k\right|^t_s\right)=\tilde{\E}\left(\left.\tilde{\gamma}^k\tilde{M}^k\right|^t_s\right)\\
    0&=\E\left(\left.\gamma^k\left(M^k\right)^2\right|^t_s-\gamma^k\int^t_s\mathbf{Q}(\n u^k(r))^2dr\right)=\tilde{\E}\left(\left.\tilde{\gamma}^k\left(\tilde{M}^k\right)^2\right|^t_s-\tilde{\gamma}^k\int^t_s\mathbf{Q}(\n \tilde{u}^k(r))^2dr\right).
  \end{aligned}
\end{equation}
The Burkholder-Davis-Gundy inequality for $W^m$ yields the uniform bound
\[\tilde{\E}\sup_{r\in[0,s]}|\tilde{W}^m(r)|^3=\E\sup_{r\in[0,s]}|W^m(r)|^3\leq Ct^{3/2}.\]
Now we can pass to the limit in (\ref{28}) by the Vitali convergence theorem and get 
\begin{equation}\label{Levy characterization theorem}
  \tilde{\E}\left(\left.\tilde{\gamma}\tilde{W}\right|^t_s\right)=0,\quad\tilde{\E}\left(\left.\tilde{\gamma}\tilde{W}^2\right|^t_s\right)=t-s.
\end{equation}
The Kolmogorov continuity theorem implies that $\tilde{W}$ has a version with continuous paths in $\R$, which we again denote by $\tilde{W}$.
Since (\ref{Levy characterization theorem}) holds for all $\gamma$, we conclude that $\tilde{W}$ is a square-integrable $(\tilde{\mathcal{F}}_t)$-martingale with its quadratic variation 
\[
  \langle\!\langle \tilde{W}\rangle\!\rangle_t=t.
\] 
Since $\tilde{W}$ is continuous, by L\'{e}vy characterization theorem, we conclude that $\tilde{W}$ is a $(\tilde{\mathcal{F}}_t)$-Wiener process on $\R$.

Now, we prove that $\tilde{M}^k(t)\rightharpoonup \tilde{M}(t)\; \p$-a.s. in $H^2(\T^n,d\mu)$. The Burkholder-Davis-Gundy inequality for $M^k$ implies for some $q>2$, and all $e\in H^2(\T^n,d\mu)$
\begin{align*}
  \tilde{\E}\sup_{s\in [0,t]}|\tilde{M}^k(s)|^q&=\E\sup_{s\in [0,t]}|M^k(s)|^q\\
  &\leq C\E\left[\int^t_0\|\mathbf{Q}(\n u^k(s))\|^2_{H^1(\T^n,d\mu)}\|e\|^2_{H^1(\T^n,d\mu)}ds\right]^{q/2}\\
  &\leq C\E\left[\int^t_0\left(1+\|\n u^k(s)\|^2_{H^1(\T^n,d\mu)}\right)\|e\|^2_{H^1(\T^n,d\mu)}ds\right]^{q/2}\\
  &\leq C_t\left(1+\E\int^t_0\|\n u^k(s)\|_{H^1(\T^n,d\mu)}ds+\sup_{s\in [0,t]}\|\n u^k(s)\|^q_{L^2(\T^n,d\mu)}\right)\|e\|^q_{H^1(\T^n,d\mu)}\\
  &\leq C\|e\|^q_{H^1(\T^n,d\mu)},
\end{align*}
uniformly in $k$. Thus $\tilde{M}^k(s)$ is uniformly integrable in $L^2(\Omega)$ with respect to $s\in [0,t]$ and all $e\in H^2(\T^n,d\mu)$. 
Since $\tilde{W}$ is a $(\tilde{\mathcal{F}}_t)$-Wiener process on $\R$, by Vitali convergence theorem, we can pass the limit in equations (\ref{29}) and infer
\begin{equation}
  \begin{aligned}
    0&=\tilde{E}\left(\tilde{\gamma}\left.\tilde{M}\right|^t_s\right).\\
    0&=\tilde{\E}\left(\left.\tilde{\gamma}\left(\tilde{M}\right)^2\right|^t_s-\tilde{\gamma}\int^t_s\mathbf{Q}(\n \tilde{u}(r))^2dr\right).
  \end{aligned}
\end{equation}
We conclude that $\tilde{M}$ is a $(H^2(\T^n,d\mu),w)$-valued continuous square-integrable $(\tilde{\mathcal{F}}_t)$-martingale and 
\[
  \left\langle\!\!\!\left\langle \tilde{M}-\int^{\cdot}_t\left\langle \mathbf{Q}(\n \tilde{u}(s))d\tilde{W}(s),e\right\rangle_{H^1(\T^n,d\mu)}\right\rangle\!\!\!\right\rangle=0.
\]
Thus,
\[
  \tilde{M}(t)-\tilde{M}(s)=\int^t_s\left\langle \mathbf{Q}(\n \tilde{u}(r))d\tilde{W}(r),e\right\rangle_{H^1(\T^n,d\mu)}dr\quad \p\text{-a.e.}
\]
for all $0\leq s\leq t<\infty$ and $e\in H^2(\T^n,d\mu)$. Furthermore we have 
\[
  \Lambda \leftharpoonup \p\circ u^k(0)^{-1}=\tilde{\p}\circ\tilde{u}^k(0)^{-1}\rightharpoonup \tilde{\p}\circ \tilde{u}(0)^{-1}. 
\]
Because of the uniform bound of $(u^{\epsilon})$ in $L^2(\Omega,L^2(0,T;H^2(\T^n,d\mu)))$ for all $T>0$, we know that the limit $\tilde{u}$ is a martingale solution of (\ref{SfMCF}) accroding to Definition \ref{definition}.

\end{proof}

\section{Small Noise Perturbation Limit of Stochastic $f$-Mean Curvature Flow}\label{sec:7}
In this section, we consider the limiting behavior of the following stochastic $f$-mean curvature flow as the noise intensity $\lambda\to0$.
\begin{equation}
  d\phi_t(x)=\no(x,t)\left(-H_f(x,t)dt+\circ \lambda dW(t)\right),\quad x\in M.
\end{equation}
As in the previous chapters, we assume $M$ is the graph of a function defined on the $n$-dimensional torus $\T^n$, i.e.,
\[
  M(t,\omega)=\left\{(x,u(x,t,\omega))\in\R^{n+1}:x\in\T^n,t\geq 0,\omega\in\Omega\right\}.
\]
In this case, we obtain the following SPDE for the height function $u$:
\begin{equation}\label{lambda equation}
  du=\left(\Q\Div_f\V+\xi u\right)dt+\lambda\Q\circ dW.
\end{equation}
When $\lambda=0$, this equation becomes the deterministic PDE for $u$:
\begin{equation}
    \partial_tu=\Q\Div_f\V+\xi u.
\end{equation}

Similar to the previous section, we also consider the following small perturbation viscous equation:
\begin{equation}\label{lambda viscous equation}
  du=\left(\epsilon\Delta_fu+\Q\Div_f\V+\xi u\right)dt+\lambda\Q\circ dW.
\end{equation}
Analogous to the proof of Theorem \ref{existence of viscous martingale solution}, we use Theorem \ref{existence} to obtain the existence of a martingale solution for the small perturbation viscous equation \eqref{lambda viscous equation}:
\begin{theorem}
  Let $\epsilon>0, q>2$, and let $\Lambda$ be a Borel probability measure on $H^1(\T^n,d\mu)$ satisfying
  \[
    \int_{H^1(\T^n,d\mu)}\|z\|^2_{H^1}d\Lambda(z)<\infty,
  \]
  and
  \[
    \int_{H^1(\T^n,d\mu)}\|\n z\|_{L^2,d\mu}^q<\infty.
  \]
  Under the conditions:
  \[D^2f\geq 0,\quad \xi\geq 0,\quad \lambda^2<2.\]
  there exists a martingale solution $u$ on $I=[0,\infty)$ with initial value $\Lambda$ for the small perturbation viscous equation (\ref{lambda viscous equation}).
\end{theorem}

Similar to Proposition \ref{uniform estimate of energy}, we obtain the following uniform energy estimate with respect to $\epsilon\geq0$:
\begin{proposition}\label{lambda uniform estimate}
  Let $\epsilon\geq 0$, and let $u$ be a solution to the small perturbation viscous equation (\ref{lambda viscous equation}). Under the pinching condition: there exists a constant $\delta>0$ such that
  \[D^2f\geq \delta\xi\mathrm{Id}> 0,\]
  and $\lambda^2<2$, the energy $\int_{\mathbb{T}^n}|\nabla u|^2d\mu$ is a supermartingale.

  More explicitly,
  \begin{align*}
        \mathbb{E}\|\nabla u\|^2_{L^2(\mathbb{T}^n,d\mu)}&+2\epsilon\mathbb{E}\int^t_0\int_{\mathbb{T}^n}|D^2u|^2d\mu ds+\left(1-\frac{\lambda^2}{2}\right)\mathbb{E}\int^t_0\int_{\mathbb{T}^n}Q^2|\operatorname{div}_f v|^2d\mu ds\\
        &+\E\int^t_0\int_{\T^n}\n u\cdot\left(\left(1+\frac{3\alpha}{2}\right)D^2f-2\alpha\xi\mathrm{Id}\right)\n ud\mu ds\\
        +\mathbb{E}\int^t_0&\int_{\mathbb{T}^n}\left(\left(1+\frac{\lambda^2}{2}\right)|D^2u|^2-\lambda^2|D^2uv|^2-\left(1-\frac{\lambda^2}{2}\right)|v\cdot D^2u v|^2\right)d\mu ds\\
        \leq \mathbb{E}\|\nabla u_0\|&_{L^2(\mathbb{T}^n,d\mu)} .
      \end{align*}
  where $\alpha$ satisfies $0<\dfrac{4\alpha}{2+3\alpha}\leq \delta$.
  Furthermore, for $q\in [1,2)$ there exists a universal constant $C>0$ such that
  \[
    \mathbb{E}\sup_{t\in I}\|\nabla u\|^{2q}_{L^2(\mathbb{T}^n,d\mu)}\leq\left(2+\frac{2C^2}{2q-q^2}\right)\mathbb{E}\|\nabla u_0\|^{2q}_{L^2(\mathbb{T}^n,d\mu)}.
  \]
  If additionally $\mathbb{P}$-$\operatorname{esssup}\|\nabla u\|_{L^{\infty}(I;L^{\infty}(\mathbb{T}^n))}=L<\infty$, then for any $t \in I$ it holds that
  \begin{align*}
    \mathbb{E}\|\nabla u(t)\|_{L^2\left(\mathbb{T}^n,d\mu\right)}^2 & +\frac{2+\lambda^2+4 L^2}{2\left(1+L^2\right)^2} \mathbb{E} \int_0^t \int_{\mathbb{T}^n}\left|D^2 u(s)\right|^2 \mathrm{~d} \mu ds \leq \mathbb{E}\left\|\nabla u_0\right\|_{L^2\left(\mathbb{T}^n,d\mu\right)}^2.
    \end{align*}
\end{proposition}
With the above uniform estimates for $\epsilon\geq0$, we can use the method from before to obtain the existence of a martingale solution for \eqref{lambda equation}:
\begin{theorem}[Existence of Martingale Solution]\label{lambda existence}
  Let $\Lambda$ be a Borel probability measure on $H^1(\T^n,d\mu)$ with finite second moment, and suppose there exists a constant $L>0$ such that the support condition holds:
  \[
    \mathrm{supp}\Lambda\subset \{z\in H^1(\T^n,d\mu): \|\nabla z\|_{L^\infty}\leq L\}.
  \]
  Then for $I=[0,\infty)$, under the pinching condition: there exists $\delta>0$ such that
  \[D^2f\geq \delta\xi\mathrm{Id}> 0\]
  and $\lambda ^2<2$, there exists a martingale solution $u^{\lambda}$ of (\ref{lambda equation}) with initial value $\Lambda$. All such solutions satisfy $$D^2u^{\lambda}\in L^2(\Omega; L^2(0,\infty;L^2(\T^n,d\mu)).$$ and
  \[
    \|\n u^{\lambda}\|_{L^{\infty}(0,\infty;L^{\infty})}\leq L\quad \mathbb{P}\text{-a.s.}.
  \]
\end{theorem}
Note that the estimates in Proposition \ref{lambda uniform estimate} are uniform with respect to both $\epsilon\geq 0$ and $\lambda ^2<2$. Therefore, we obtain the following uniform estimates for $u^{\lambda}$ with respect to $\lambda^2<2$:
\begin{proposition}
   Let $\lambda^2<2$, and let $u^{\lambda}$ be the martingale solution from Theorem \ref{lambda existence}. Under the pinching condition: there exists a constant $\delta>0$ such that
   \[ D^2f\geq \delta\xi \mathrm{Id}>0.
   \] it holds that
   \begin{align*}
        \mathbb{E}\|\nabla u^{\lambda}\|^2_{L^2(\mathbb{T}^2,d\mu)}
        &+\left(1-\frac{\lambda^2}{2}\right)\mathbb{E}\int^t_0\int_{\mathbb{T}^n}Q^2|\operatorname{div}_f v|^2d\mu ds\\
        &+\E\int^t_0\int_{\T^n}\n u^{\lambda}\cdot\left(\left(1+\frac{3\alpha}{2}\right)D^2f-2\alpha\xi\mathrm{Id}\right)\n u^{\lambda}d\mu ds\\
        +\mathbb{E}\int^t_0\int_{\mathbb{T}^n}&\left(\left(1+\frac{\lambda^2}{2}\right)|D^2u^{\lambda}|^2-\lambda^2|D^2u^{\lambda}v|^2-\left(1-\frac{\lambda^2}{2}\right)|v\cdot D^2u^{\lambda} v|^2\right)d\mu ds\\
        \leq \mathbb{E}\|\nabla u_0\|&_{L^2(\mathbb{T}^n,d\mu)}.
      \end{align*}
      where $\alpha$ satisfies $0<\dfrac{4\alpha}{2+3\alpha}\leq \delta$. For $q\in [1,2)$ there exists a universal constant $C>0$ such that
      \[
    \mathbb{E}\sup_{t\in I}\|\nabla u^{\lambda}\|^{2q}_{L^2(\mathbb{T}^n,d\mu)}\leq\left(2+\frac{2C^2}{2q-q^2}\right)\mathbb{E}\|\nabla u_0\|^{2q}_{L^2(\mathbb{T}^n,d\mu)}.
  \]
  and
  \begin{align*}
    \mathbb{E}\|\nabla u(t)\|_{L^2\left(\mathbb{T}^n,d\mu\right)}^2 & +\frac{2+\lambda^2+4 L^2}{2\left(1+L^2\right)^2} \mathbb{E} \int_0^t \int_{\mathbb{T}^n}\left|D^2 u(s)\right|^2 \mathrm{~d} \mu ds \leq \mathbb{E}\left\|\nabla u_0\right\|_{L^2\left(\mathbb{T}^n,d\mu\right)}^2   .
    \end{align*}
\end{proposition}
With the above uniform estimates for $\lambda^2<2$, we can use the proof method from the previous section to show that the distributions of $(u^{\lambda})_{\lambda^2<2}$ are tight in $\mathcal{X}_u$. By the Jakubowski-Skorohkod representation theorem, $(u^{\lambda})_{\lambda^2<2}$ is relatively compact. For any sequence $\lambda _k\to 0$, by relative compactness, there exists a subsequence (still denoted by $(\lambda_k)$) such that $u^{\lambda_k}$ converges almost surely. It is straightforward to verify that this limit satisfies the following deterministic $f$-mean curvature flow equation:
\[
\partial_tu=\Q\Div_f\V+\xi u.
\]
If this equation has a unique solution $u$, then we have
\[
u^{\lambda }\to u \text{\quad a.s. in\quad }\mathcal{X}_u\text{\quad as\quad}\lambda \to 0.
\]
Thus we obtain the following theorem:
\begin{theorem}
   Let $\lambda^2<2$, and let $u^{\lambda}$ be the martingale solution of the stochastic $f$-mean curvature flow
   \begin{equation*}
  du=\left(\Q\Div_f\V+\xi u\right)dt+\lambda\Q\circ dW.
\end{equation*} from Theorem \ref{lambda existence}. Under the pinching condition: there exists a constant $\delta>0$ such that \[
D^2f\geq \delta\xi\mathrm{Id}>0.
\]
If the deterministic $f$-mean curvature flow equation:
\[
\partial_tu=\Q\Div_f\V+\xi u,
\]
has a unique solution $u$, then we have
\[
u^{\lambda }\to u\in \mathcal{X}_u \text{\quad a.s. \quad as\quad}\lambda \to 0.
\]
\end{theorem}
\begin{remark}
    If $f$ in our stochastic $f$-mean curvature flow is constant, then it reduces to the stochastic mean curvature flow
    \[
    du=\left(\Q\mathrm{div}\dfrac{\n u}{\Q}\right)dt+\lambda\Q\circ dW.
    \] The results of the above theorem also hold. In particular, \cite{DS23} provides criteria for the uniqueness of solutions to the mean curvature flow equation
    \[
    \partial_t u=\Q\mathrm{div}\dfrac{\n u}{\Q}.
    \]
\end{remark}


\bibliographystyle{plain}
\bibliography{export}

\begin{bibdiv}
\begin{biblist}

\bib{BM14}{incollection}{
      author={Borisenko, Aleksander},
      author={Rovenski, Vladimir},
       title={Gaussian mean curvature flow for submanifolds in space forms},
        date={2014},
   booktitle={Geometry and its applications},
      series={Springer Proc. Math. Stat.},
      volume={72},
   publisher={Springer, Cham},
       pages={39\ndash 50},
         url={https://doi.org/10.1007/978-3-319-04675-4_2},
      review={\MR{3213506}},
}

\bib{BM10}{article}{
      author={Borisenko, Alexander~A.},
      author={Miquel, Vicente},
       title={Gaussian mean curvature flow},
        date={2010},
        ISSN={1424-3199,1424-3202},
     journal={J. Evol. Equ.},
      volume={10},
      number={2},
       pages={413\ndash 423},
         url={https://doi.org/10.1007/s00028-010-0054-2},
      review={\MR{2643802}},
}

\bib{Giga}{article}{
      author={Chen, Yun~Gang},
      author={Giga, Yoshikazu},
      author={Goto, Shun'ichi},
       title={Uniqueness and existence of viscosity solutions of generalized mean curvature flow equations},
        date={1991},
        ISSN={0022-040X,1945-743X},
     journal={J. Differential Geom.},
      volume={33},
      number={3},
       pages={749\ndash 786},
         url={http://projecteuclid.org/euclid.jdg/1214446564},
      review={\MR{1100211}},
}

\bib{Prévôt2007}{book}{
      author={Claudia, Prévôt},
      author={Michael, Röckner},
       title={A concise course on stochastic partial differential equations},
   publisher={Springer Berlin Heidelberg},
     address={Berlin, Heidelberg},
        date={2007},
        ISBN={978-3-540-70781-3},
         url={https://doi.org/10.1007/978-3-540-70781-3_4},
}

\bib{Colding2010}{article}{
      author={Colding, Tobias~Holck},
      author={II, William P~Minicozzi},
      author={Pedersen, Erik~Kjaer},
       title={Mean curvature flow},
        date={2010},
     journal={Bull. Amer. Math. Soc. (N.S.)},
      volume={52},
}

\bib{daprato_zabczyk_2014}{book}{
      author={Da~Prato, Giuseppe},
      author={Zabczyk, Jerzy},
       title={Stochastic equations in infinite dimensions},
     edition={2},
      series={Encyclopedia of Mathematics and its Applications},
   publisher={Cambridge University Press},
        date={2014},
}

\bib{Dabrock2020}{thesis}{
      author={Dabrock, Nils},
       title={Stochastic mean curvature flow},
        type={PhD thesis, Fakult\"{a}t Mathematik, TU Dortmund},
        date={2020},
         url={http://dx.doi.org/10.17877/DE290R-21118},
}

\bib{Dabrock2021}{article}{
      author={Dabrock, Nils},
      author={Hofmanov\'{a}, Martina},
      author={R\"{o}ger, Matthias},
       title={Existence of martingale solutions and large-time behavior for a stochastic mean curvature flow of graphs},
        date={2021},
        ISSN={0178-8051,1432-2064},
     journal={Probab. Theory Related Fields},
      volume={179},
      number={1-2},
       pages={407\ndash 449},
         url={https://doi.org/10.1007/s00440-020-01012-6},
      review={\MR{4221662}},
}

\bib{DS23}{article}{
      author={Daskalopoulos, Panagiota},
      author={Saez, Mariel},
       title={Uniqueness of entire graphs evolving by mean curvature flow},
        date={2023},
        ISSN={0075-4102,1435-5345},
     journal={J. Reine Angew. Math.},
      volume={796},
       pages={201\ndash 227},
         url={https://doi.org/10.1515/crelle-2022-0085},
      review={\MR{4554469}},
}

\bib{Dirr2001}{article}{
      author={Dirr, Nicolas},
      author={Luckhaus, Stephan},
      author={Novaga, Matteo},
       title={A stochastic selection principle in case of fattening for curvature flow},
        date={2001},
        ISSN={0944-2669,1432-0835},
     journal={Calc. Var. Partial Differential Equations},
      volume={13},
      number={4},
       pages={405\ndash 425},
         url={https://doi.org/10.1007/s005260100080},
      review={\MR{1867935}},
}

\bib{Huisken1989}{article}{
      author={Ecker, Klaus},
      author={Huisken, Gerhard},
       title={Mean curvature evolution of entire graphs},
        date={1989},
        ISSN={0003-486X,1939-8980},
     journal={Ann. of Math. (2)},
      volume={130},
      number={3},
       pages={453\ndash 471},
         url={https://doi.org/10.2307/1971452},
      review={\MR{1025164}},
}

\bib{Renesse2012}{article}{
      author={Es-Sarhir, Abdelhadi},
      author={von Renesse, Max-K.},
       title={Ergodicity of stochastic curve shortening flow in the plane},
        date={2012},
        ISSN={0036-1410,1095-7154},
     journal={SIAM J. Math. Anal.},
      volume={44},
      number={1},
       pages={224\ndash 244},
         url={https://doi.org/10.1137/100798235},
      review={\MR{2888287}},
}

\bib{Gawarecki2011}{book}{
      author={Gawarecki, Leszek},
      author={Mandrekar, Vidyadhar},
       title={Stochastic differential equations in infinite dimensions: with applications to stochastic partial differential equations},
   publisher={Springer Berlin Heidelberg},
     address={Berlin, Heidelberg},
        date={2011},
        ISBN={978-3-642-16194-0},
         url={https://doi.org/10.1007/978-3-642-16194-0_4},
}

\bib{Hofmanova2017}{article}{
      author={Hofmanov\'{a}, Martina},
      author={R\"{o}ger, Matthias},
      author={von Renesse, Max},
       title={Weak solutions for a stochastic mean curvature flow of two-dimensional graphs},
        date={2017},
        ISSN={0178-8051,1432-2064},
     journal={Probab. Theory Related Fields},
      volume={168},
      number={1-2},
       pages={373\ndash 408},
         url={https://doi.org/10.1007/s00440-016-0713-5},
      review={\MR{3651056}},
}

\bib{Huisken1984}{article}{
      author={Huisken, Gerhard},
       title={Flow by mean curvature of convex surfaces into spheres},
        date={1984},
        ISSN={0022-040X,1945-743X},
     journal={J. Differential Geom.},
      volume={20},
      number={1},
       pages={237\ndash 266},
         url={http://projecteuclid.org/euclid.jdg/1214438998},
      review={\MR{772132}},
}

\bib{Jakubowski1997}{article}{
      author={Jakubowski, A.},
       title={The almost sure {S}korokhod representation for subsequences in nonmetric spaces},
        date={1997},
        ISSN={0040-361X},
     journal={Teor. Veroyatnost. i Primenen.},
      volume={42},
      number={1},
       pages={209\ndash 216},
         url={https://doi.org/10.1137/S0040585X97976052},
      review={\MR{1453342}},
}

\bib{JL06}{article}{
      author={Jian, Huai~Yu},
      author={Liu, Yan~Nan},
       title={Ginzburg-{L}andau vortex and mean curvature flow with external force field},
        date={2006},
        ISSN={1439-8516,1439-7617},
     journal={Acta Math. Sin. (Engl. Ser.)},
      volume={22},
      number={6},
       pages={1831\ndash 1842},
         url={https://doi.org/10.1007/s10114-005-0698-y},
      review={\MR{2262444}},
}

\bib{JL08}{article}{
      author={Jian, Huai-Yu},
      author={Liu, Yannan},
       title={Long-time existence of mean curvature flow with external force fields},
        date={2008},
        ISSN={0030-8730,1945-5844},
     journal={Pacific J. Math.},
      volume={234},
      number={2},
       pages={311\ndash 325},
         url={https://doi.org/10.2140/pjm.2008.234.311},
      review={\MR{2373451}},
}

\bib{kawasaki1982}{article}{
      author={Kawasaki, Kyozi},
      author={Ohta, Takao},
       title={Kinetic {Drumhead} {Model} of {Interface}. {I}},
        date={1982-01},
        ISSN={0033-068X},
     journal={Progress of Theoretical Physics},
      volume={67},
      number={1},
       pages={147\ndash 163},
         url={https://doi.org/10.1143/PTP.67.147},
        note={\_eprint: https://academic.oup.com/ptp/article-pdf/67/1/147/5438937/67-1-147.pdf},
}

\bib{LW2015}{article}{
      author={Li, Haizhong},
      author={Wei, Yong},
       title={{$f$}-minimal surface and manifold with positive {$m$}-{B}akry-\'{E}mery {R}icci curvature},
        date={2015},
        ISSN={1050-6926,1559-002X},
     journal={J. Geom. Anal.},
      volume={25},
      number={1},
       pages={421\ndash 435},
         url={https://doi.org/10.1007/s12220-013-9434-5},
      review={\MR{3299288}},
}

\bib{XD-Li}{article}{
      author={Li, Xiang-Dong},
       title={Riesz transforms on forms and {$L^p$}-{H}odge decomposition on complete {R}iemannian manifolds},
        date={2010},
        ISSN={0213-2230,2235-0616},
     journal={Rev. Mat. Iberoam.},
      volume={26},
      number={2},
       pages={481\ndash 528},
         url={https://doi.org/10.4171/RMI/607},
      review={\MR{2677005}},
}

\bib{Jian2007}{article}{
      author={Liu, Yan-nan},
      author={Jian, Huai-yu},
       title={Evolution of hypersurfaces by the mean curvature minus an external force field},
        date={2007},
        ISSN={1006-9283,1862-2763},
     journal={Sci. China Ser. A},
      volume={50},
      number={2},
       pages={231\ndash 239},
         url={https://doi.org/10.1007/s11425-007-2077-x},
      review={\MR{2306089}},
}

\bib{Liu2012}{article}{
      author={Liu, Yannan},
       title={Evolution of hypersurfaces by powers of mean curvature minus an external force field},
        date={2012},
        ISSN={1673-3452,1673-3576},
     journal={Front. Math. China},
      volume={7},
      number={4},
       pages={717\ndash 723},
         url={https://doi.org/10.1007/s11464-012-0218-1},
      review={\MR{2944568}},
}

\bib{Jian2009}{article}{
      author={Liu, Yannan},
      author={Jian, Huaiyu},
       title={Evolution of spacelike hypersurfaces by mean curvature minus external force field in {M}inkowski space},
        date={2009},
        ISSN={1536-1365,2169-0375},
     journal={Adv. Nonlinear Stud.},
      volume={9},
      number={3},
       pages={513\ndash 522},
         url={https://doi.org/10.1515/ans-2009-0305},
      review={\MR{2536952}},
}

\bib{meira_space_2020}{article}{
      author={Meira, Adson},
      author={Gon\c{c}alves, Rosivaldo~Antonio},
       title={On the space of {$f$}-minimal surfaces with bounded {$f$}-index in weighted smooth metric spaces},
        date={2020},
        ISSN={0025-2611,1432-1785},
     journal={Manuscripta Math.},
      volume={162},
      number={3-4},
       pages={559\ndash 563},
         url={https://doi.org/10.1007/s00229-019-01144-7},
      review={\MR{4109500}},
}

\bib{Viot}{inproceedings}{
      author={Metivier, Michel},
      author={Viot, Michel},
       title={On weak solutions of stochastic partial differential equations},
        date={1988},
   booktitle={Stochastic analysis},
      editor={M{\'e}tivier, Michel},
      editor={Watanabe, Shinzo},
   publisher={Springer Berlin Heidelberg},
     address={Berlin, Heidelberg},
       pages={139\ndash 150},
}

\bib{Pardoux2021}{book}{
      author={Pardoux, {\'E}tienne},
       title={Spdes as infinite-dimensional sdes},
   publisher={Springer International Publishing},
     address={Cham},
        date={2021},
        ISBN={978-3-030-89003-2},
         url={https://doi.org/10.1007/978-3-030-89003-2_2},
}

\bib{Seidler1993}{article}{
      author={Seidler, Jan},
       title={Da prato-zabczyk's maximal inequality revisited. i.},
    language={eng},
        date={1993},
     journal={Mathematica Bohemica},
      volume={118},
      number={1},
       pages={67\ndash 106},
         url={http://eudml.org/doc/29167},
}

\bib{Simon}{article}{
      author={Simon, Jacques},
       title={Compact sets in the space {$L^p(0,T;B)$}},
        date={1987},
        ISSN={0003-4622},
     journal={Ann. Mat. Pura Appl. (4)},
      volume={146},
       pages={65\ndash 96},
         url={https://doi.org/10.1007/BF01762360},
      review={\MR{916688}},
}

\bib{Smoczyk01}{article}{
      author={Smoczyk, Knut},
       title={A relation between mean curvature flow solitons and minimal submanifolds},
        date={2001},
        ISSN={0025-584X,1522-2616},
     journal={Math. Nachr.},
      volume={229},
       pages={175\ndash 186},
         url={https://doi.org/10.1002/1522-2616(200109)229:1<175::AID-MANA175>3.3.CO;2-8},
      review={\MR{1855161}},
}

\bib{Souganidis2004}{article}{
      author={Souganidis, P.~E.},
      author={Yip, N.~K.},
       title={Uniqueness of motion by mean curvature perturbed by stochastic noise},
        date={2004},
        ISSN={0294-1449,1873-1430},
     journal={Ann. Inst. H. Poincar\'{e} C Anal. Non Lin\'{e}aire},
      volume={21},
      number={1},
       pages={1\ndash 23},
         url={https://doi.org/10.1016/S0294-1449(03)00029-5},
      review={\MR{2037245}},
}

\bib{Wei09}{article}{
      author={Wei, Guofang},
      author={Wylie, Will},
       title={Comparison geometry for the {B}akry-{E}mery {R}icci tensor},
        date={2009},
        ISSN={0022-040X,1945-743X},
     journal={J. Differential Geom.},
      volume={83},
      number={2},
       pages={377\ndash 405},
         url={https://doi.org/10.4310/jdg/1261495336},
      review={\MR{2577473}},
}

\bib{Yip1998}{article}{
      author={Yip, Nung~Kwan},
       title={Stochastic motion by mean curvature},
        date={1998},
        ISSN={0003-9527},
     journal={Arch. Rational Mech. Anal.},
      volume={144},
      number={4},
       pages={313\ndash 355},
         url={https://doi.org/10.1007/s002050050120},
      review={\MR{1656479}},
}

\end{biblist}
\end{bibdiv}

\end{document}